\documentclass[12pt,draftclsnofoot,onecolumn]{IEEEtran}

\usepackage{psfrag}
\usepackage{booktabs}
\usepackage{epsfig}
\usepackage{latexsym}
\usepackage{amsmath}
\usepackage{amssymb}
\usepackage{amsfonts}
\usepackage{bm}
\usepackage{colortbl}
\usepackage{soul}
\usepackage{hhline}
\usepackage{cite}

\newtheorem{definition}{Definition}
\newtheorem{theorem}{Theorem}
\newtheorem{proposition}[theorem]{Proposition}
\newtheorem{corollary}[theorem]{Corollary}
\newtheorem{lemma}[theorem]{Lemma}
\newtheorem{assumption}{Assumption}

\newcommand {\otoprule }{\midrule[\heavyrulewidth]}

\renewcommand{\Re}{{\mathbb R}}
\def\QED{~\rule[-1pt]{5pt}{5pt}\par\medskip}
\newenvironment{proof_ref}[1]{\noindent {\it Proof of #1\/:}}{ \hfill{\QED}}
\newenvironment{proof_text}{\noindent {\it Proof:}}{\hfill {\QED}}

\begin{document}
\allowdisplaybreaks

\title{Stochastic optimization on continuous domains
with finite-time guarantees
by Markov chain Monte Carlo methods}

\author{A. Lecchini-Visintini$^*$\footnote{$^*$Department of Engineering,
  University of Leicester, {\tt alv1@leicester.ac.uk}.},
  J. Lygeros$^\dag$\footnote{$^\dag$Automatic Control
  Laboratory, ETH Zurich, {\tt lygeros@control.ee.ethz.ch}.} and
  J. Maciejowski$^\ddag$\footnote{$^\ddag$Department of Engineering, University of
  Cambridge {\tt jmm@eng.cam.ac.uk}.}}

\date{\today}
\maketitle
\vspace{-1cm}
\begin{abstract}
We introduce bounds on the finite-time performance of
Markov chain Monte Carlo algorithms in  approaching the global
solution of stochastic optimization problems
over continuous domains.\linebreak
A comparison with other state-of-the-art methods having finite-time guarantees
for solving stochastic programming problems is  included.
\end{abstract}

\section{Introduction}\label{sec:Intro}

 In principle, any optimization problem on a finite domain can be
solved by an exhaustive search. However, this is often beyond
computational capacity: the optimization domain of the traveling
salesman problem with $100$ cities contains more than $10^{155}$
possible tours \cite{Applegate-et-al-06}. An efficient algorithm to solve the traveling
salesman and many similar problems has not yet been found and such
problems remain solvable only in principle
\cite{Achlioptas-et-al-05}. Statistical mechanics has inspired
widely used methods for finding good approximate solutions in hard
discrete optimization problems which defy efficient exact solutions
\cite{Kirkpatrick-et-al-83,Bonomi-Lutton-84,Fu-Anderson-86,Mezard-et-al-02}.
Here a key idea has been that of simulated annealing
\cite{Kirkpatrick-et-al-83}: a random search based on the
Metropolis-Hastings algorithm, such that the distribution of the
elements of the domain visited during the search converges to an
equilibrium distribution concentrated around the global optimizers.
Convergence and finite-time performance of simulated annealing on
finite domains have been evaluated \mbox{e.g.~in \cite{Laarhoven-Aarts-87,Mitra-et-al-86,Hajek-88,Hannig-et-al-06}.}

On continuous domains, most popular optimization methods perform a
local gradient-based search and in general converge to local
optimizers, with the notable exception of convex optimization problems where
convergence to the unique global optimizer occurs \cite{Boyd-Vandenberghe-03}.
Simulated annealing performs a global
search and can be easily implemented on continuous domains
using the general family of Markov chain Monte Carlo (MCMC) methods \cite{Robert-Casella-04}.
Hence it can be considered a powerful complement to local methods. In this
paper,  we introduce  for the first time rigorous guarantees on the
finite-time performance of simulated annealing on continuous
domains. We will show that it is possible to derive  MCMC algorithms
to implement simulated annealing
which can find an approximate solution to the problem of
optimizing a function of continuous variables, within a specified
tolerance and with an arbitrarily high level of
confidence after a known finite number
of steps. Rigorous guarantees on the finite-time performance of
simulated annealing in the optimization of functions of continuous
variables have never been obtained before; the only results
available state that simulated annealing converges to a global
optimizer as the number of steps grows to infinity,
e.g.~\cite{Haario-Saksman-91,Gelfand-Mitter-91,Tsallis-Stariolo-96,Locatelli-00,Andrieu-et-al-01},
asymptotic convergence rates have been \mbox{obtained in \cite{Locatelli-01,Rubenthaler-et-al-09}.}

The background of our work is twofold. On the one hand, our definition of
``approximate domain optimizer", introduced in Section \ref{sec:Opt} as an approximate solution to a global
optimization problem, is inspired by the definition of ``probably approximate near minimum" introduced
by Vidyasagar in \cite{Vidyasagar-01} for global
optimization based on the concept of finite-time learning with known accuracy
and confidence of statistical learning \mbox{theory \cite{Vapnik-95,Vidyasagar-03}.}
In the control field the work of Vidyasagar \cite{Vidyasagar-01,Vidyasagar-03} has been seminal
in the development of the so-called randomized approach.
Inspired by statistical learning theory,
this approach is characterized by the construction of algorithms
which make use of independent sampling in order to find probabilistic approximate solutions
to difficult control system design applications
see e.g.~\cite{Tempo-et-al-05,Calafiore-Campi-06,Alamo-et-al-07} and the references therein.
In our work, the definition of approximate domain optimizer will be essential in establishing
rigorous guarantees on the finite-time performance of simulated annealing.
On the other hand, we show that our rigorous finite-time guarantees
can be achieved  by the wider class of algorithms based on
Markov chain Monte Carlo sampling.
Hence, we ground our results on the theory of convergence, with
quantitative bounds on the distance to the target distribution, of
the Metropolis-Hastings algorithm and MCMC methods \cite{Meyn-Tweedie-93,Rosenthal-95,Mengersen-Tweedie-96,Roberts-Rosenthal-04}.
In addition, we demonstrate how, under some quite weak regularity conditions,
our definition of approximate domain optimizer
can be related to the standard notion of approximate optimization
considered in the stochastic programming literature  \cite{Shapiro-Nemirovski-05,Shapiro-08,Nesterov-05,Nesterov-Vial-08}.
This link provides theoretical support for the use of
simulated annealing and MCMC optimization algorithms, which
have been proposed, for example, in~\cite{Muller-99,Doucet-et-al-02,Muller-et-al-04},
for solving stochastic programming problems.
In this paper, beyond the presentation of some simple illustrative examples,
we will not develop any  ready-to-use optimization algorithm.
The Metropolis-Hastings algorithm and the
general family of MCMC  methods have many degrees of freedom. The
choice and comparison of specific algorithms goes beyond the scope
of the paper.

The paper is organized as follows. In Section \ref{sec:Opt} we introduce
the definition of approximate domain optimizer and establish a direct
relationship between the approximate domain optimizer and the standard notion of
approximate optimizer adopted in the stochastic programming literature.
In Section \ref{sec:MCMC} we first recall the reasons why existing results on the convergence
of simulated  annealing on continuous domain do not provide finite-time guarantees.
Then we state the main results of the paper
and we discuss their consequences.
In Section \ref{sec:PCA} we illustrate the convergence of MCMC algorithms.
In Section \ref{sec:NE} we present a simple illustrative numerical example.
In Section \ref{sec:Disc} we compare the MCMC approach with other state-of-the-art
methods for solving stochastic programming problems with finite-time
performance bounds.
In  Section \ref{sec:end} we state our findings and conclude the paper.
The Appendix contains all technical proofs.
Some of the results of this paper were included in preliminary
conference contributions \cite{Lecchini-et-al-08,Lecchini-et-al-08-II}.

\section{Approximate optimizers}\label{sec:Opt}

Consider an optimization criterion $U:\Theta \rightarrow \Re$,
with $\Theta\subseteq \Re^n$, and let
\begin{equation}
U^*:=\sup_{\theta \in \Theta}U(\theta). \label{eq:sup0}
\end{equation}
The following will be a standing assumption for all our results.\smallskip
\begin{assumption}\label{assu:1}
$\Theta$ has finite Lebesgue measure. $U$ is well defined point-wise,
measurable, and bounded between 0 and 1 (i.e.~$U(\theta)\in[0,\, 1]\,\,\forall\theta\in\Theta$).
\end{assumption}\smallskip
In general, any bounded criterion can be scaled to take values in $[0,\, 1]$.
Given, for example, $U'(\theta) \in [\underline{U},\, \overline{U}]$
we can consider the optimization of the modified function
$$
U(\theta) =
\frac{U'(\theta)-\underline{U}}{\overline{U}-\underline{U}}\, ,
$$
which takes values in $[0, 1]$ for all $\theta \in \Theta$.
(In this case, we need to multiply the value imprecision, $\epsilon$ below,
by $(\overline{U}-\underline{U})$ to obtain
its corresponding value in the scale of the original criterion $U'$.)

For some results another assumption will be needed.\smallskip
\begin{assumption}\label{assu:2}
$\Theta$ is compact.
$U$ is Lipschitz continuous.
\end{assumption}\smallskip
We use $L$ to denote the Lipschitz constant of $U$,
i.e.~$\forall \theta_1, \theta_2 \in \Theta, \; |U(\theta_1)-U(\theta_2)| \leq L \|\theta_1-\theta_2\|$.
Assumption~\ref{assu:2} implies the existence of a global optimizer,
i.e.~under Assumption~\ref{assu:2}, we have $ \Theta^*:=\{\theta \in \Theta \; | \; U(\theta)=U^*\} \neq \emptyset.$

If, given an element $\theta$ in $\Theta$, the
value $U(\theta)$ can be computed directly, we say that $U$ is a
deterministic criterion. In this case the optimization
problem~\eqref{eq:sup0} is a standard, in general non-linear,
non-smooth, programming problem. Examples of such a deterministic
optimization criterion are, among many possible others, the design
criterion in a robust control design problem \cite{Vidyasagar-01} and
the energy landscape in protein structure prediction~\cite{Wales-03}.
In problems involving random
variables, the value $U(\theta)$ can be
the expected value of some function $g:\Theta\times X \rightarrow\Re$
which depends on both the optimization variable
$\theta$, and on some random variable $\bm{x}$ with probability
distribution $P_{\bm{x}}(\cdot;\theta)$ which may itself depend on
$\theta$, i.e.
\begin{equation}\label{eq:exval}
U(\theta) = \int g(x,\theta)P_{\bm{x}}(dx;\theta)\, .
\end{equation}
In such problems it is usually not possible to compute
$U(\theta)$ directly. In stochastic optimization
\cite{Shapiro-Nemirovski-05,Shapiro-08,Nesterov-05,Nesterov-Vial-00,Muller-99,Doucet-et-al-02,Muller-et-al-04},
it is typically assumed that one can  obtain independent samples of $\bm{x}$ for a
given $\theta$, hence obtain sample values of $g(\bm{x},\theta)$,
and thus construct a Monte Carlo estimate of $U(\theta)$.
In some application it might not be possible or efficient to obtain independent
samples of $\bm{x}$. In this case one has to resort to other Monte Carlo
strategies to approximate $U(\theta)$ such as, for example, importance sampling \cite{Robert-Casella-04}.
The Bayesian experimental design of clinical trials is an important application area
where expected-value criteria arise \cite{Spiegelhalter-et-al-04}.
We investigate the optimization of expected-value criteria
motivated by problems of aircraft routing \cite{Lecchini-et-al-06}
and parameter identification for genetic networks~\cite{KCL07FOSBE}.
In the particular case that $P_{\bm{x}}(dx;\theta)$ does not depend on
$\theta$, the `$\inf$' counterpart of problem  (\ref{eq:sup0})
is called ``empirical risk minimization'', and is studied extensively in statistical learning
theory \cite{Vapnik-95,Vidyasagar-03}.
Conditions on $g$ and $P_{\bm{x}}$ to ensure that
$U$ is Lipschitz continuous (for Assumption \ref{assu:2}) can be found in~\cite[pag.~189-190]{Shapiro-08}.
The results reported here
apply in the same way to the optimization of both deterministic and
expected-value criteria.

 We introduce two different definitions of
approximate solution to the optimization \mbox{problem (\ref{eq:sup0}).}
The first is the definition of approximate domain optimizer.
It will be essential in establishing finite-time guarantees on the
performance of MCMC methods.
\smallskip
\begin{definition}
\label{def:DomOpt}
Let $\epsilon\geq0$ and  $\alpha\in [0,\,1]$ be given numbers.
Then $\theta$ is an approximate domain optimizer of $U$ with value imprecision
$\epsilon$ and residual domain $\alpha$ if
\begin{equation}\label{eq:DomOpt}
\lambda(\{\theta'\in\Theta :  U(\theta') > U(\theta) +\epsilon \}) \leq \alpha\, \lambda(\Theta)
\end{equation}
where $\lambda$ denotes the Lebesgue measure.
\end{definition}\smallskip
That is, the function $U$ takes values strictly greater than $U(\theta) + \epsilon$ only
on a subset of values of $\theta$ no larger than an $\alpha$ portion
of the optimization domain. The smaller $\epsilon$  and
$\alpha$ are, the better is the approximation of a true global
optimizer. If  both $\alpha$ and $\epsilon$ are equal to zero then
$U(\theta)$ coincides with the  essential supremum of $U$ \cite{EssentialSupremum}.
We will use
$$
\Theta(\epsilon,\alpha):=\{\theta \in \Theta \; | \;
\lambda(\{\theta' \in \Theta \; | \; U(\theta')>U(\theta)+\epsilon\})
\leq \alpha \lambda(\Theta)\}
$$
to denote the set of approximate domain optimizers with value
imprecision $\epsilon$ and residual domain $\alpha$.
The intuition that our notion of approximate domain
optimizer can be used to obtain formal guarantees on the finite-time
performance of optimization methods based on a stochastic search of
the domain is already apparent in the work of Vidyasagar.
Vidyasagar \cite{Vidyasagar-03,Vidyasagar-01} introduces the
similar definition of ``probably approximate near minimum"
and obtains rigorous finite-time guarantees in the optimization of
expected value criteria based on uniform independent sampling of the
domain.
The method of Vidyasagar has had considerable success in solving difficult
control system design applications
\cite{Vidyasagar-01,Tempo-et-al-05}. Its  appeal stems from its
rigorous finite-time guarantees which exist without the need for any
particular  assumption on the optimization criterion.

The following is a more common notion of approximate optimizer.\smallskip
\begin{definition}\label{def:ValOpt}
Let $\epsilon\geq0$ be a given number. Then $\theta$ is an
 an approximate value optimizer of $U$ with imprecision $\epsilon$ if
$U(\theta') \leq U(\theta)+\epsilon$ for all $\theta' \in \Theta$.
\end{definition}\smallskip
This notion is commonly used in the stochastic programming literature
\cite{Shapiro-Nemirovski-05,Shapiro-08,Nesterov-05,Nesterov-Vial-00} and provides a direct
bound on $U^*$: $\theta\in \Theta$ is an approximate
value optimizer with imprecision $\epsilon>0$ if and only if
$U^*\leq U(\theta)+\epsilon$. We will use
$$
\Theta^*(\epsilon):=\{\theta \in \Theta \; | \;
\forall \theta' \in \Theta, \; U(\theta') \leq U(\theta)+\epsilon\}
$$
to denote the set of approximate value optimizers with imprecision $\epsilon$.

It is easy to see that for all $\epsilon$
if $\Theta^*\neq\emptyset$ then
$\Theta^* \subseteq \Theta^*(\epsilon)$.
Notice that $\Theta^*(\epsilon)$ does not coincide
with $\Theta(\epsilon,\alpha)$.
In fact one can see that approximate value optimality is a stronger concept
than approximate domain optimality, in the following sense.
For all $\epsilon$ and all $\alpha$,
if $\Theta^*(\epsilon)\neq\emptyset$ then
$\Theta^*(\epsilon)\subseteq
\Theta(\epsilon,\alpha)$.
Conversely, given an approximate domain optimizer it is in general not
possible to draw any conclusions about the approximate value
optimizers. For example, for any $\alpha$ the function $U:[0, 1]
\rightarrow [0, 1]$ with
\[
U(\theta)=\left\{\begin{array}{ll}
1 & \mbox{ if } \theta \in [0, \alpha)\\
0 & \mbox{ if } \theta \in [\alpha, 1]
\end{array}\right.
\]
has the property that $\Theta(\epsilon,\alpha)=\Theta$ for all
$\epsilon>0$. Therefore, given $\theta \in
\Theta(\epsilon,\alpha)$ it is impossible to draw any conclusions
about $U^*$; the only possible bound is $U^* \leq U(\theta)+1$
which, given that $U(\theta) \in [0, 1]$, is meaningless.
A relation between domain and value approximate optimality can,
however, be established under Assumption~\ref{assu:2}. \smallskip
\begin{theorem}\label{thm:1}
Let Assumption~\ref{assu:2} hold.  Let $\theta$ be an approximate domain optimizer
with value imprecision $\epsilon$ and residual domain $\alpha$.
Then, $\theta$ is also an approximate value optimizer with
imprecision
$$
\epsilon + \frac{L}{\sqrt{\pi}}
\left[\frac{n}{2}\Gamma\left(\frac{n}{2}\right)\right]^{\frac{1}{n}} [\alpha\lambda(\Theta)]^{\frac{1}{n}}
$$
where $\Gamma$ denotes the gamma function.
\end{theorem}\smallskip
Theorem \ref{thm:1} shows that
\begin{equation}\label{eq:th1-txt}
\theta \in \Theta(\epsilon,\alpha) \Rightarrow
U^* \leq U(\theta) + \epsilon
+\frac{L}{\sqrt{\pi}}
\left[\frac{n}{2}\Gamma\left(\frac{n}{2}\right)\right]^{\frac{1}{n}}
[\alpha\lambda(\Theta)]^{\frac{1}{n}}\ .
\end{equation}
The result allows us to select the value of $\alpha$ in such a way that
an approximate domain optimizer with value imprecision $\epsilon$
and residual domain $\alpha$ is also an
approximate value optimizer with imprecision $2\epsilon$.
To do this, we need to select $\alpha$ so that
$\frac{L}{\sqrt{\pi}}
\left[\frac{n}{2}\Gamma\left(\frac{n}{2}\right)\right]^{\frac{1}{n}}
[\alpha\lambda(\Theta)]^{\frac{1}{n}} \leq \epsilon$ hence
\begin{equation}\label{eq:alphabound}
\alpha \leq
\frac{\left[\frac{\epsilon \sqrt{\pi}}{L}\right]^n}
{\lambda(\Theta)\left[\frac{n}{2}\Gamma\left(\frac{n}{2}\right)\right]}\, .
\end{equation}
To illustrate the above inequality consider the case where the domain
$\Theta$ is contained in an $n$-dimensional ball of radius $R$.
Notice that under Assumption \ref{assu:2} the existence of such an $R$ is guaranteed.
In this case
$ \lambda(\Theta) = \frac{2\pi^{\frac{n}{2}}}{n\Gamma(\frac{n}{2})}R^n\,.$
Therefore (\ref{eq:alphabound}) becomes
\begin{equation}
\label{eq:alphaboundR}
\alpha \leq \left(\frac{1}{L}\frac{\epsilon}{R}\right)^n\ .
\end{equation}
Note that, as $n$ increases, $\alpha$ has to decrease to zero rapidly
to ensure the required imprecision of the approximate
value optimizer. In this case, $\alpha$ needs to decrease to zero as $\epsilon^n$.

\section{Optimization with MCMC: finite time guarantees}\label{sec:MCMC}

In simulated annealing, a random search based on the
Metropolis-Hastings algorithm is carried out,
such that the distribution of the
elements of the domain visited during the search converges to an
equilibrium distribution concentrated around the global optimizers.

Here we adopt equilibrium distributions
defined by densities proportional to $[U(\theta)+\delta]^J$, where $J$ and
$\delta$ are strictly positive parameters.
We use
\begin{equation}\label{eq:equi}
\pi(d\theta;J,\delta)\propto[U(\theta)+\delta]^J\lambda(d\theta)
\end{equation}
to denote this equilibrium distribution.
The presence of $\delta$ is a technical condition required in the proof
of our main result and will be discussed
later on in this section.
In our setting, the so-called `zero-temperature' distribution
is the limiting distribution $\pi(\,\cdot\,\,;J,\delta)$
for $J\rightarrow \infty$ denoted by $\pi_\infty$.
It can be shown that under some technical conditions, $\pi_\infty$ is a
uniform distribution on the set $\Theta^*$ of the global maximizers
of $U$ \cite{Hwang-80}.

\begin{figure}\centering
\includegraphics[width=0.8\columnwidth]{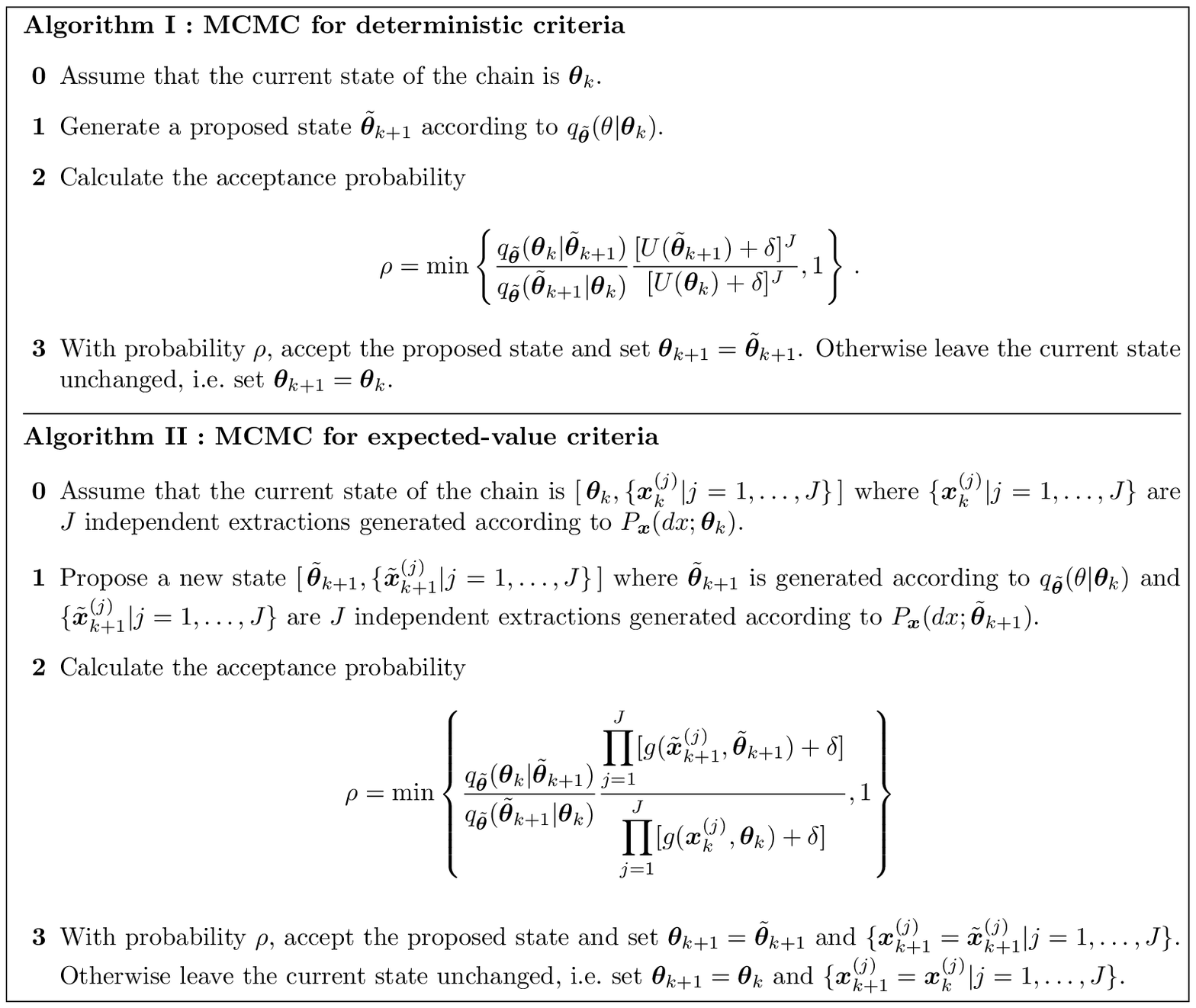}
\caption{The basic iterations of the Metropolis-Hastings algorithm with
equilibrium distributions $\pi(\cdot;J,\delta)$
for the maximization of deterministic and expected-value criteria.
In both algorithms, $q_{\tilde{\bm{\theta}}}(\cdot|\bm{\theta}_k)$
is the density of the `proposal distribution'.
}\label{fig:algo}
\end{figure}

 In Fig.~\ref{fig:algo}, we illustrate two algorithms
which implement  Markov transition kernels
with equilibrium distributions $\pi(\,\cdot\,\,;J,\delta)$.
Algorithm I is the `classical' Metropolis-Hastings algorithm
for the case in which $U$ is a deterministic criterion.
Algorithm II is a suitably modified version of the Metropolis-Hastings
algorithm  for the case in which $U$
is an expected-value criterion in the form of (\ref{eq:exval}).
This latter algorithm was devised by M\"{u}ller \cite{Muller-99, Muller-et-al-04}
and Doucet et al.~\cite{Doucet-et-al-02}.

 In the simulated annealing scheme, one would simulate an inhomogeneous chain
in which the Markov transition kernel at the $k$-th step of the chain
has equilibrium distribution $\pi(\,\cdot\,\,;J_k,\delta)$ where
 $\{J_k\}_{k=1,2,\dots}$ is a suitably chosen `cooling schedule',
i.e.~a non-decreasing sequence of values for the exponent $J$.
The rationale of simulated annealing is as follows: if the
temperature is kept constant, say $J_k=J$, then the distribution of
the state of the chain, say  $P_{\bm{\theta}_k}$, tends to the equilibrium
distribution $\pi(\,\cdot\,\,;J,\delta)$; if $J\rightarrow\infty$ then the
equilibrium distribution $\pi(\,\cdot\,\,;J,\delta)$  tends to the zero-temperature
distribution $\pi_\infty$; as a result, if the cooling schedule
$J_k$ tends to $\infty$, one obtains that the distribution of the state
of the chain  $P_{\bm{\theta}_k}$ tends to $\pi_\infty$
\cite{Haario-Saksman-91,Gelfand-Mitter-91,Tsallis-Stariolo-96,Locatelli-00,Andrieu-et-al-01}.

 The difficulty which must be overcome in order to obtain finite step
results on simulated annealing algorithms
on a continuous domain is
that usually, in an optimization problem defined over continuous
variables, the set of global optimizers $\Theta^*$ has zero Lebesgue measure
(e.g.~a set of isolated points).
Notice that this is not the case
for a finite domain, where the set of global optimizers is of
non-null measure with respect to the reference counting measure
\cite{Laarhoven-Aarts-87,Mitra-et-al-86,Hajek-88,Hannig-et-al-06}.
It is instructive to look at the issue in terms of the rate of
convergence to the target zero-temperature distribution.
On a continuous domain, the standard distance between
two distributions, say $\mu_1$ and $\mu_2$, is the total variation distance
$\|\mu_1 - \mu_2\|_{\mbox{\tiny TV}} = \sup_{A\in\mathfrak{B}(\Theta)} | \mu_1(A) - \mu_2(A)|$.
If the  set of global optimizers $\Theta^*$ has
 zero Lebesgue measure, then the target
zero-temperature distribution $\pi_\infty$ ends up being a
mixture of probability masses on $\Theta^*$.
On the other hand, the distribution of the state of the chain
$P_{\bm{\theta}_k}$ is absolutely continuous with respect to the
Lebesgue measure (i.e.~\mbox{$\lambda(A)=0\Rightarrow
P_{\bm{\theta}_k}(A)=0$}) by construction for any finite $k$.
Hence, if $\Theta^*$ has zero Lebesgue measure then it has
zero measure also according to $P_{\bm{\theta}_k}$.
The set $\Theta^*$ has however measure 1  according to $\pi_\infty$.
The distance $\|P_{\bm{\theta}_k} - \pi_\infty\|_{\mbox{\tiny TV}}$ is
then constantly $1$.
In general, on a continuous domain, although the distribution of the state of the chain
$P_{\bm{\theta}_k}$  converges asymptotically  to $\pi_\infty$,
it is not possible to introduce a sensible
distance between $P_{\bm{\theta}_k}$ and $\pi_\infty$
and a rate of convergence to the target distribution cannot even be defined
(weak convergence), see \cite[Theorem 3.3]{Haario-Saksman-91}.

Weak convergence to $\pi_\infty$  implies that,
asymptotically, $\bm{\theta}_k$ eventually hits the set of approximate value optimizers
$\Theta^*(\epsilon)$, for any $\epsilon>0$, with probability one \cite{Haario-Saksman-91,Gelfand-Mitter-91,Tsallis-Stariolo-96,Locatelli-00,Andrieu-et-al-01}.
In more recent works,  bounds on the expected number of iterations before hitting
$\Theta^*(\epsilon)$ \cite{Locatelli-01} or on
$P_{\bm{\theta}_k}(\Theta^*(\epsilon))$ \cite{Rubenthaler-et-al-09} have been obtained.
In \cite{Rubenthaler-et-al-09}, a short review of existing bounds is proposed,
and  under some technical conditions, it is proven that
for any $\epsilon>0$ there is a number $C_\epsilon$  such that
$P_{\bm{\theta}_k}(\{\theta \in \Theta \; | \; U(\theta)\leq U^* - \epsilon \})\leq C_\epsilon k^{-\frac{1}{3}}(1+\log k)$.
In general, the expressions in these bounds cannot be computed.
For example, in the bound reported here, $C_\epsilon$ is not known in advance.
Hence, existing bounds can be used to asses the asymptotic
rate of convergence  but not as stopping criteria.

 Here we show that finite-time guarantees for stochastic optimization
by MCMC methods on continuous domains can
be obtained by selecting a distribution $\pi(\,\cdot\,\,;J,\delta)$
with a finite $J$ as the target distribution in place of the
zero-temperature distribution $\pi_\infty$.
Our definition of approximate domain optimizer given
in Section \ref{sec:Opt} is essential for establishing this result.
The definition of approximate domain optimizers carries
an important property, which holds regardless of what the
criterion $U$ is: if $\epsilon$ and $\alpha$ have non-zero values
then the set of approximate global optimizers
$\Theta(\epsilon,\alpha)$ always has non-zero Lebesgue measure.
The following theorem establishes a lower bound on the measure of
the set $\Theta(\epsilon,\alpha)$ with respect
to a distribution $\pi(\cdot;J,\delta)$ with finite $J$.
It is important to stress that the result holds universally for {\it any}
optimization criterion $U$ on a bounded domain. The only  minor
requirement is that $U$ takes values in $[0,\, 1]$.\smallskip
\begin{theorem}\label{thm:bound}
Let Assumption~\ref{assu:1} hold.
Let $\Theta(\epsilon,\alpha)$ be the set of approximate
domain optimizers of $U$ with value imprecision $\epsilon$
and residual domain $\alpha$.
Let $J\geq1$ and $\delta>0$, and consider the distribution
$\pi(d\theta;J,\delta)\propto [U(\theta)+\delta]^J\lambda(d\theta)$.
Then,  for any $\alpha\in (0,\, 1]$ and $\epsilon\in [0,\, 1]$,
the following inequality holds
\begin{equation}\label{eq:confidence}
\pi(\Theta(\epsilon,\alpha);J,\delta)\geq
\frac { 1 }
{ \mbox{$\displaystyle 1 + \left[\frac{1 + \delta}{\epsilon + 1 + \delta}\right]^{\,J}
\left[\frac{1}{\alpha}\frac{1 +\delta}{\epsilon+\delta} - 1\right]
\frac{1+\delta}{\delta} $}}\,\, .
\end{equation}
\end{theorem}\smallskip
Notice that, for given
non-zero values of $\epsilon$, $\alpha$, and $\delta$
the right-hand side of (\ref{eq:confidence}) can be
made arbitrarily close to 1 by choice of $J$.
To obtain some insight on this choice it is instructive to turn the
bound of Theorem \ref{thm:bound}
around to provide a lower bound on $J$ which ensures
that $\pi(\Theta(\epsilon,\alpha);J,\delta)$ attains some desired value $\sigma$.\smallskip
\begin{corollary}\label{cor:Jbound}
Let the notation and assumptions of Theorem \ref{thm:bound} hold.
For any  $\alpha\in (0,\, 1]$, $\epsilon\in (0,\, 1]$
and $\sigma\in(0,\,1)$, if
\begin{equation}\label{eq:Jbound}
J \geq \frac{1+\epsilon+\delta}{\epsilon}
\left[\log\frac{\sigma}{1-\sigma}+
\log\frac{1}{\alpha}+2\log\frac{1+\delta}{\delta}\right]
\end{equation}
then $\pi(\Theta(\epsilon,\alpha);J,\delta) \geq \sigma$.
\end{corollary}\smallskip
The importance of the choice of a target distribution $\pi(\,\cdot\,\,;J,\delta)$
with a finite $J$ is that the distance
$\|P_{\bm{\theta}_k}  - \pi(\,\cdot\,\,;J,\delta)\|_{\mbox{\tiny TV}}$
is a meaningful quantity.
Convergence of the Metropolis-Hastings algorithm and MCMC methods in
total variation distance is a well studied problem. The theory provides
simple conditions under which one derives upper bounds on
$\|P_{\bm{\theta}_k}  - \pi(\,\cdot\,\,;J,\delta)\|_{\mbox{\tiny TV}}$ that
decrease to zero as  $k\rightarrow\infty$ \cite{Meyn-Tweedie-93,Rosenthal-95,Mengersen-Tweedie-96,Roberts-Rosenthal-04}.
 It is then appropriate to introduce the following finite-time result.\smallskip
\begin{proposition}\label{prop:confidence_k}
Let the notation and assumptions of Theorem \ref{thm:bound} hold.
Assume that $J$ respects the bound of Corollary \ref{cor:Jbound}
for given $\alpha$, $\epsilon$, $\delta$ and $\sigma$.
Let $\bm{\theta}_k$ with distribution $P_{\bm{\theta}_k}$ be
the state of the chain of an MCMC algorithm
with target distribution $\pi(\,\cdot\,\,;J,\delta)$.
Then,
$$
P_{\bm{\theta}_k}(\Theta(\epsilon,\alpha);J,\delta)\geq
\sigma - \|P_{\bm{\theta}_k} - \pi(\,\cdot\,\,;J,\delta)\|_{\mbox{\tiny TV}}\, .
$$
In other words, the statement ``$\bm{\theta}_k$ is an approximate domain optimizer of $U$ with
value imprecision  $\epsilon$ and residual domain $\alpha$" can be made
with confidence  $\sigma-\|P_{\bm{\theta}_k} - \pi(\,\cdot\,\,;J,\delta)\|_{\mbox{\tiny TV}}$.
\end{proposition}\smallskip
The proof follows directly from the definition of the
total variation distance.

If the optimization criterion is Lipschitz continuous,
Theorem \ref{thm:bound} can be used together with
Theorem~\ref{thm:1} to derive a lower bound on the measure of the
set of approximate value optimizers with a given imprecision
with respect to a distribution $\pi(\,\cdot\,\,;J,\delta)$.
An example of such a bound is the following.\smallskip
\begin{proposition}\label{prop:confidence_kR}
Let the notation and assumptions of Theorems \ref{thm:1} and \ref{thm:bound} hold.
In addition, assume that  $\Theta$ is contained in an $n$-dimensional ball of radius $R$.
Let $\bm{\theta}_k$ with distribution $P_{\bm{\theta}_k}$ be
the state of the chain of an MCMC algorithm
with target distribution $\pi(\,\cdot\,\,;J,\delta)$.
For given $\epsilon\in (0,\, 1]$ and $\sigma\in(0,\,1)$, if
\begin{equation}
\label{eq:JboundvalueR}
J \geq \frac{1+\epsilon+\delta}{\epsilon}
\left[\log\frac{\sigma}{1-\sigma}+
n\log\left(\frac{LR}{\epsilon}\right)+2\log\frac{1+\delta}{\delta}\right]
\end{equation}
then
$$
P_{\bm{\theta}_k}(\Theta^*(2\epsilon);J,\delta)\geq
\sigma - \|P_{\bm{\theta}_k} - \pi(\,\cdot\,\,;J,\delta)\|_{\mbox{\tiny TV}}\, .
$$
In other words, the statement
``$\bm{\theta}_k$ is an approximate value optimizer of $U$ with
value imprecision  $2\epsilon$"  can be made
with confidence  $\sigma-\|P_{\bm{\theta}_k} - \pi(\,\cdot\,\,;J,\delta)\|_{\mbox{\tiny TV}}$.
\end{proposition}\smallskip
The proof follows by substituting $\alpha$ with
the right-hand side of (\ref{eq:alphaboundR}) in (\ref{eq:Jbound})
and from the definition of the total variation distance.

Finally, Theorem~\ref{thm:bound} provides a criterion for selecting
the parameter $\delta$ in $\pi(\,\cdot\,\,;J,\delta)$.
For given $\epsilon$ and $\alpha$, there exists
an optimal choice of $\delta$ which minimizes the value of $J$
required to ensure $\pi(\Theta(\epsilon,\alpha);J,\delta) \geq \sigma$.
The advantage of choosing  the smallest $J$,
consistent with the required $\sigma$, is computational.
The exponent $J$ coincides with the number of
Monte Carlo simulations of random variable $\bm{x}$  which must
be done at each step in Algorithm II.
The smallest $J$ reduces also the peakedness of $\pi(\cdot;J,\delta)$.
The higher the peakedness of $\pi(\cdot;J,\delta)$ is the harder is
to design a proposal distribution which operates efficiently.
In turn, reducing the peakedness of $\pi(\cdot;J,\delta)$
will decrease the number of steps required to achieve
the desired reduction of $\|P_{\bm{\theta}_k}  - \pi(\,\cdot\,\,;J,\delta)\|_{\mbox{\tiny TV}}$.
\begin{figure}[t]\centering
\includegraphics[width=0.8\columnwidth]{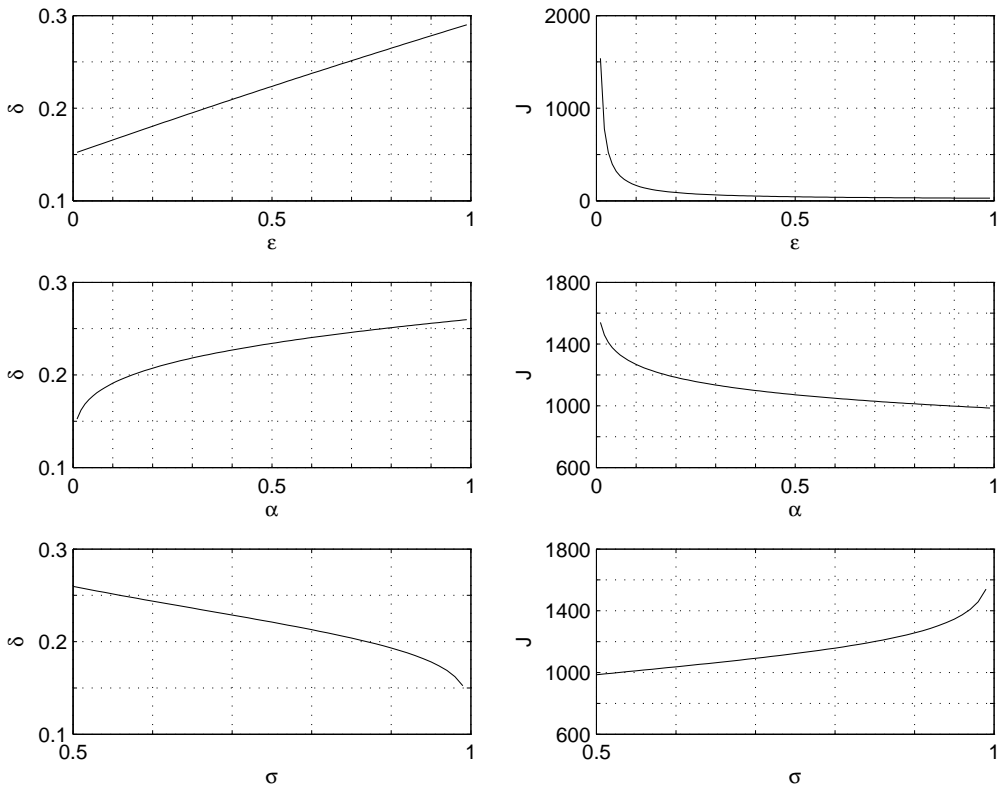}
\caption{Variation of the optimal $\delta$ and of the corresponding $J$
with respect to $\epsilon$, $\alpha$ and $\sigma$.
Two of the three parameters are kept constant for each figure,
set to $\epsilon=0.01$, $\alpha=0.01$ and $\sigma=0.99$.}\label{fig:deltaJ}
\end{figure}
The optimal choice of $\delta$ is specified by the following result. \smallskip
\begin{proposition}\label{prop:optidelta}
For fixed $\epsilon>0$, $\alpha>0$, and $\sigma\in(0.5,\,1)$, the function
\[
f(\delta) =
\frac{1+\epsilon+\delta}{\epsilon}\left[\log\frac{\sigma}{1-\sigma}+
  \log\frac{1}{\alpha}+2\log\frac{1+\delta}{\delta}\right]\, ,
\]
i.e.~the right hand side of inequality (\ref{eq:Jbound}),
is convex in $\delta$ and attains its global minimum at the
unique solution (for $\delta$) of the equation
\[
\log\frac{1+\delta}{\delta} + \log\frac{\sqrt\sigma}{\sqrt{1-\sigma}} +
\log\frac{1}{\sqrt{\alpha}} =
\frac{1+\epsilon+\delta}{\delta(1+\delta)}.
\]
\end{proposition}\smallskip
For example, if $\epsilon=0.01$, $\alpha=0.01$ and  $\sigma=0.99$,
then one obtains $\delta = 0.15$ and $J=1540$.
Plots of the value of the optimal $\delta$ and of the corresponding
value of $J$ for different values of $\epsilon$, $\alpha$ and
$\sigma$ are shown in Fig.~\ref{fig:deltaJ}.
Notice that the result of Proposition \ref{prop:optidelta} holds also
for inequality (\ref{eq:JboundvalueR})
provided that $\alpha$ in the statement of Proposition \ref{prop:optidelta}
is replaced by the right hand side of (\ref{eq:alphaboundR}).

\section{Convergence}\label{sec:PCA}

In this section we illustrate the statement of Proposition \ref{prop:confidence_k}.
We base the discussion on the simplest
available result on the convergence of MCMC methods
in total variation distance, taken from~\cite{Mengersen-Tweedie-96}.
In this case,  the proposal distribution, denoted by its density
$q_{\tilde{\bm{\theta}}}(\theta|\bm{\theta}_k)$
in Algorithms I and II,
is independent of the current state $\bm{\theta}_k$.\smallskip
\begin{theorem}[\cite{Mengersen-Tweedie-96}]\label{thm:MT}
Let $P_{\bm{\theta}_k}$ be the distribution of the state of the chain
in the Metropolis-Hastings algorithm with an independent proposal distribution.
Let $\pi$ denote the target distribution.
Let $p$ and $q$ denote respectively the density of $\pi$ and
the density of the proposal distribution and assume that
$p(\theta)>0,\,\forall\theta\in\Theta$ and $q(\theta)>0,\,\forall\theta\in\Theta$.
If there exists $M$ such that
$p(\theta)\leq M q(\theta)\,,\,\,\,\forall \theta\in\Theta$,
then
\begin{equation}\label{eq:MT}
\left\|\pi-P_{\bm{\theta}_k}\right\|_{\mbox{\tiny TV}} \leq \left(1-\frac{1}{M}\right)^k.
\end{equation}
\end{theorem}
{\it Proof:} See \cite[Theorem 2.1]{Mengersen-Tweedie-96}, or~\cite[Theorem 7.8]{Robert-Casella-04}.\smallskip\\
Here, we chose $q_{\tilde{\bm{\theta}}}$ as the uniform distribution over $\Theta$.
Sampling using an independent uniform proposal distribution is a na\"{\i}ve strategy in
an MCMC approach and cannot be expected to perform efficiently \cite{Robert-Casella-04}.
However, it allows us to present some simple illustrative examples where convergence
bounds can be derived with a few basic steps.

In some cases the na\"{\i}ve strategy  can produce approximate domain optimizers
very efficiently. One such case occurs under the assumption that the
optimization criterion $U(\theta)$ has a ``flat top'',  i.e. the set
of global optimizers $\Theta^*$ has non-zero Lebesgue measure. The
same assumption has been used in \cite[Theorem 4.2]{Haario-Saksman-91}
to obtain the strong convergence of simulated annealing on a continuous domain.
In this case, the application of Theorem \ref{thm:MT} provides the following result.\smallskip
\begin{proposition}\label{prop:MTJb}
Let the notation and assumptions of Proposition \ref{prop:confidence_k} hold.
In particular, assume that $\bm{\theta}_k$ is the state of the chain
of the Metropolis-Hastings algorithm
with independent uniform proposal distribution.
In addition, given $\rho\in(0,\,1)$, let $\sigma = (1+\gamma)\rho$
for some  $\gamma\in(0,\, \frac{1-\rho}{\rho})$.
Let $\Theta^*$ be the set of global optimizer of $U$ and assume
that $\lambda(\Theta^*)\geq\beta\lambda(\Theta)$ for some $\beta\in(0, 1)$.
If
\begin{equation}\label{eq:MTJb}
k \geq \frac{\log{\gamma\rho}}{\log(1-\beta)}
\end{equation}
then  $P_{\bm{\theta}_k}(\Theta(\epsilon,\alpha);J,\delta)\geq \rho$.
\end{proposition}\smallskip
In (\ref{eq:MTJb}), it is convenient to choose  $\gamma\approx\frac{1-\rho}{\rho}$.
Hence, the number of iterations grows approximately
as  $-\log(1-\rho)=\log(\frac{1}{1-\rho})$ and $-\frac{1}{\log(1-\beta)}$
and is independent of $\epsilon$ and $\alpha$.
In Algorithm II the total number of required
samples of $\bm{x}$ is given by the number of iterations multiplied by $J$.
In this case, it can be shown that a nearly optimal choice is $\gamma=\frac{1}{2}\frac{1-\rho}{\rho}$.
Hence, using (\ref{eq:Jbound}) for the case of approximate domain optimization,
we obtain that the required samples of $\bm{x}$ grow as
$\frac{1}{\epsilon}$, $\log\frac{1}{\alpha}$, and approximately as $(\log\frac{1}{1-\rho})^2$.
Instead,  using (\ref{eq:JboundvalueR}) for the case of approximate value optimization,
we obtain that the required samples of $\bm{x}$ grow as
$\frac{1}{\epsilon}\log\frac{1}{\epsilon}$, $(\log\frac{1}{1-\rho})^2$,
\mbox{$\log LR$ and $n$.}

If the `flat top' condition is not met it can be easily
seen that the use of a uniform proposal
distribution can lead to an exponential number of iterations.
The problem is the implicit dependence of the
convergence rate on the exponent $J$.
In the general case, by applying Theorem \ref{thm:MT} we obtain the following result.\smallskip
\begin{proposition}\label{prop:MTJ}
Let the notation and assumptions of Proposition \ref{prop:confidence_k} hold.
In particular, assume that $\bm{\theta}_k$ is the state of the chain
of the Metropolis-Hastings algorithm
with independent uniform proposal distribution.
In addition, given $\rho\in(0,\,1)$, let $\sigma=(1+\gamma)\rho$
for some  $\gamma\in(0,\, \frac{1-\rho}{\rho})$.
If $k \geq \left(\frac{1+\delta}{\delta}\right)^J\log\left(\frac{1}{\gamma\rho} \right) $
 or, equivalently,
\begin{equation}\label{eq:MTJ}
k \geq \left[\frac{(1+\gamma)\rho}{1-(1+\gamma)\rho}\frac{1}{\alpha}
  \left(\frac{1+\delta}{\delta}\right)^2
  \right]^{\frac{1+\epsilon+\delta}{\epsilon}\log\left(\frac{1+\delta}{\delta}\right)}
  \log\frac{1}{\gamma\rho}
\end{equation}
then $P_{\bm{\theta}_k}(\Theta(\epsilon,\alpha);J,\delta)\geq \rho$.
\end{proposition}\smallskip
Hence, the number of iterations turns out to be exponential in $\frac{1}{\epsilon}$.
In Algorithm II, the total number of required extractions of $\bm{x}$ grows like $Jk$,
which is also exponential in $\frac{1}{\epsilon}$.
Therefore, using Theorem \ref{thm:MT} for Algorithms I and II with
 $q_{\tilde{\bm{\theta}}}$ as the independent
uniform proposal distribution, the only general bounds that we
can guarantee are exponential.

\section{Numerical example}\label{sec:NE}

\begin{figure}[t]
\centerline{
  \includegraphics[width=0.5\hsize]{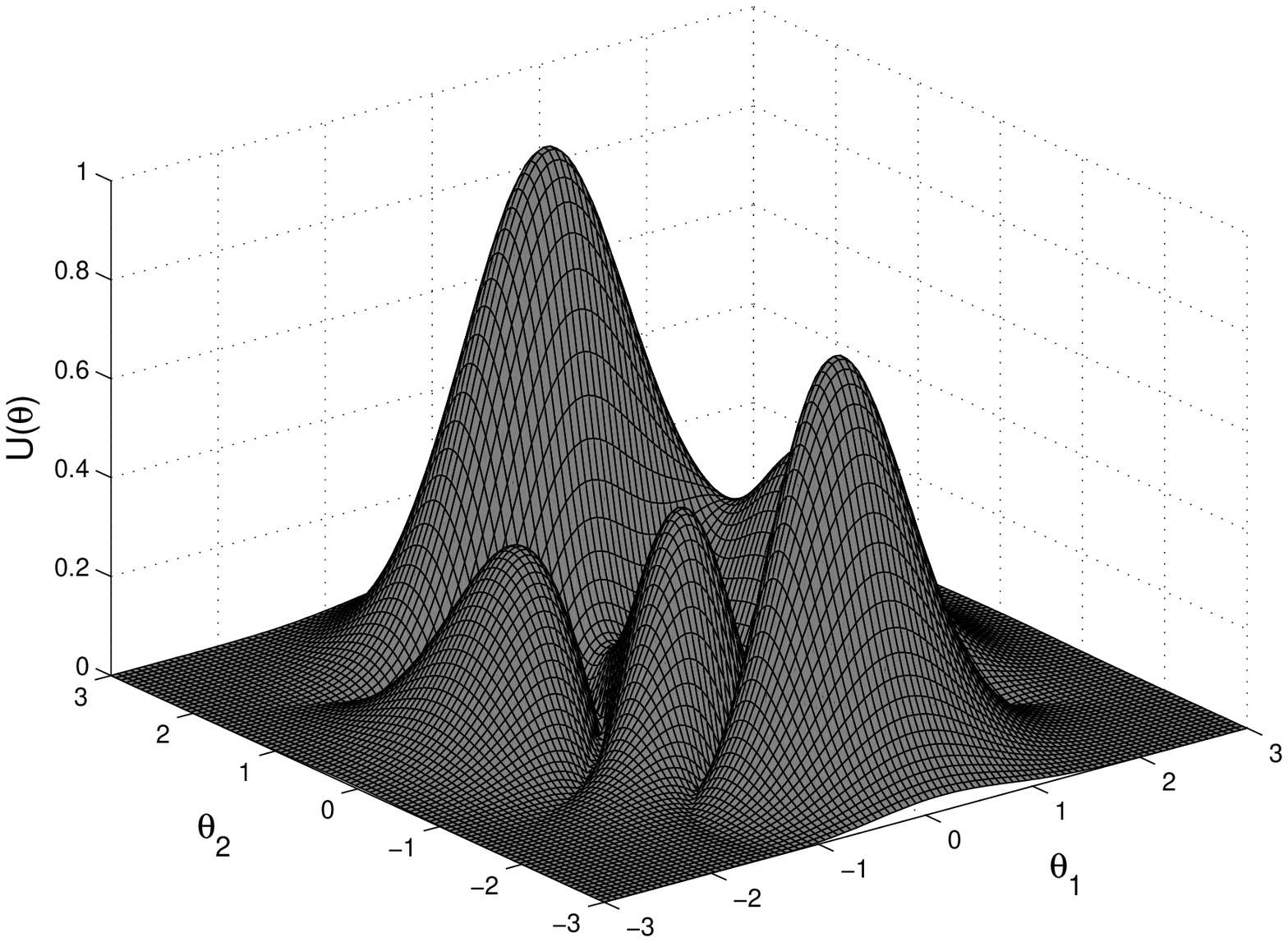}
  \includegraphics[width=0.5\hsize]{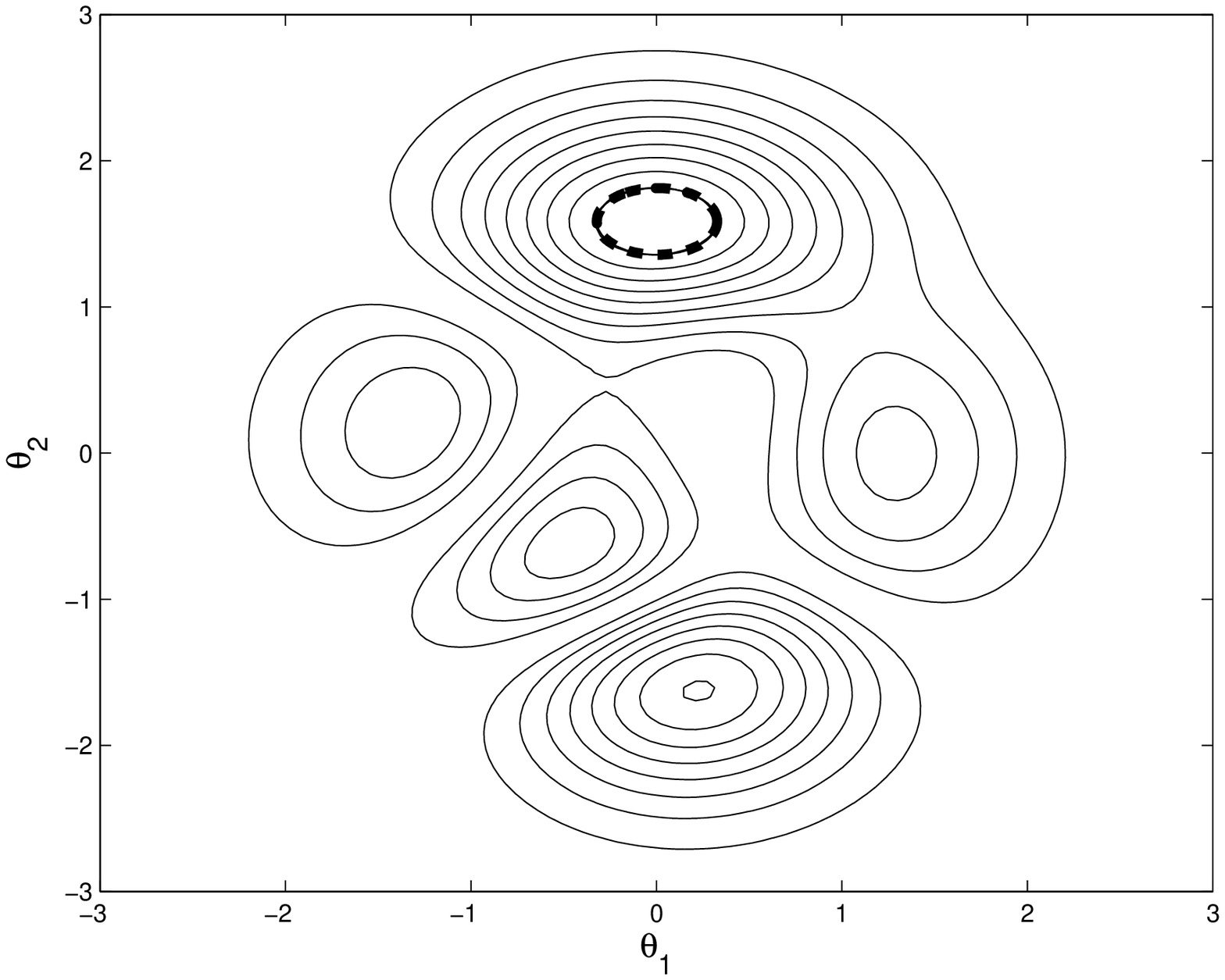}
}
\caption{Function $U(\theta)$ (left panel) and its level sets (right panel).
The $0.9$ level set is highlighted as a dashed ellipse.}
\label{fig:fun}
\end{figure}

To demonstrate some of the bounds derived in this work we
apply the proposed method to a simple example.
Let $\theta \in \Theta = [-3,3] \times [-3,3]$ and consider the function
\[
V(\theta)=3(1-\theta_1)^2e^{-\theta_1^2-(\theta_2+1)^2}
-10(\frac{\theta_1}{5}-\theta_1^3-\theta_2^5)e^{-\theta_1^2-\theta_2^2}
-\frac{1}{3}e^{-(\theta_1+1)^2-\theta_2^2}
\]
(the Matlab function {\tt peaks}).
We define the function $U: \Theta \rightarrow [0,1]$ by
\[
U(\theta)=\frac{|V(\theta)|}{\max_{\theta' \in \Theta}|V(\theta')|}.
\]
The scaling factor $\max_{\theta' \in \Theta}|V(\theta')|=8.1062$ and
a Lipschitz constant of $U(\theta)$, $L=1.725$, were computed numerically
using a grid on $\Theta$. The function $U$ and its level sets are shown
in Fig.~\ref{fig:fun}. The $0.9$
level set, which coincides with $\Theta^*(0.1)$, is highlighted in the figure.

To obtain a stochastic programming problem multiplicative noise was
added using the function
\[
g(\bm{x},\theta)=(1+\bm{x})U(\theta)
\]
where $\bm{x}$ is normally distributed with mean $0$ and variance $0.25$. It
is easy to see that the expected value of $g(\bm{x},\theta)$ is indeed equal
to $U(\theta)$. One can think of $g(\bm{x},\theta)$ as an imperfect, unbiased
measurement device of $U(\theta)$: We can only collect information about $U$
through noise corrupted samples generated by $g$. Notice that the noise
intensity is higher in areas where $U(\theta)$ is large, making it more
difficult to use the samples to pinpoint the maxima of $U$.

The MCMC Algorithm II of Fig.~\ref{fig:algo} was applied to this function.
The design parameter $\delta=0.1$ and an
independent uniform proposal distribution $q$ were used throughout.

\begin{figure}[t!]
\centerline{
  \includegraphics[width=0.5\hsize]{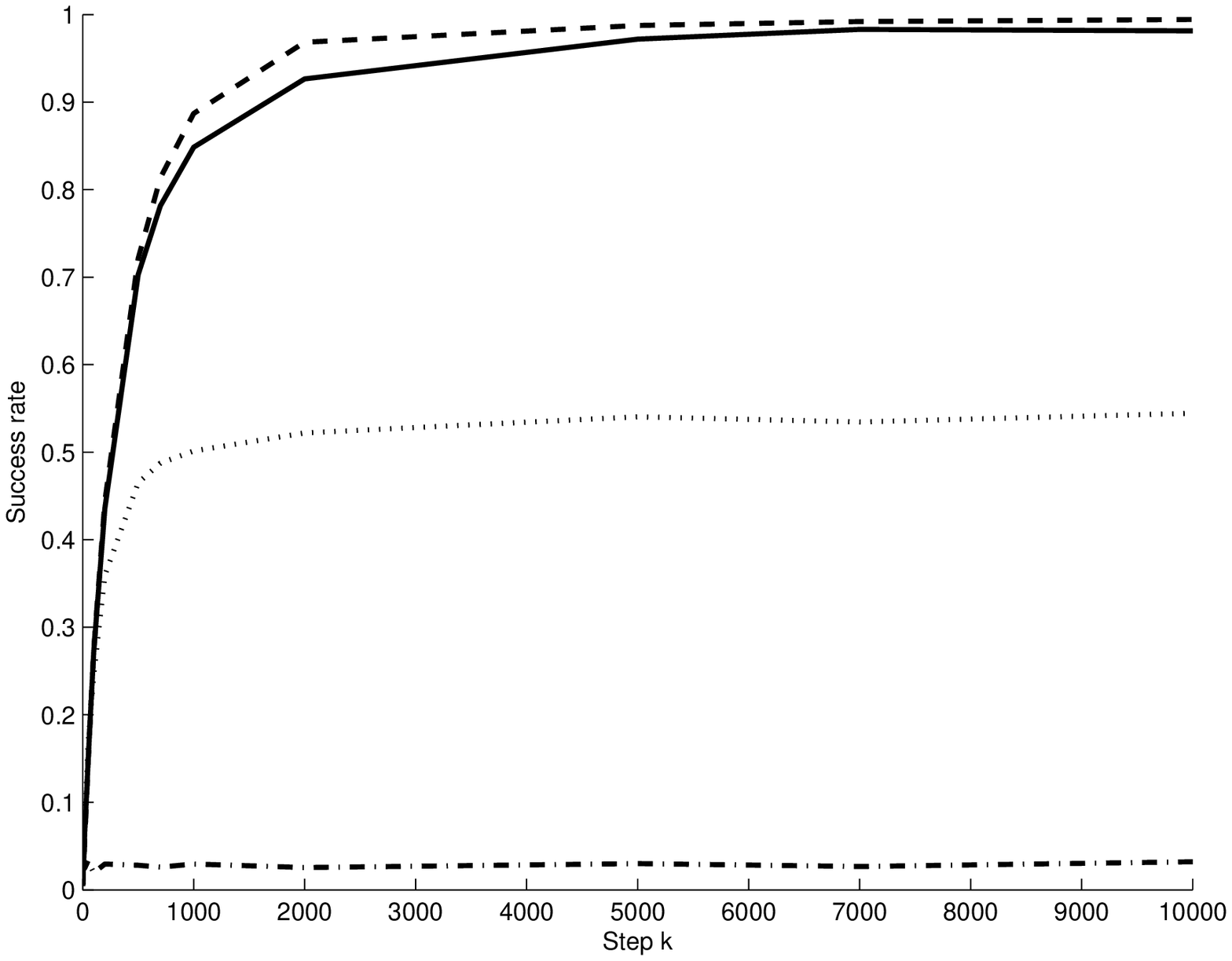}
  \includegraphics[width=0.5\hsize]{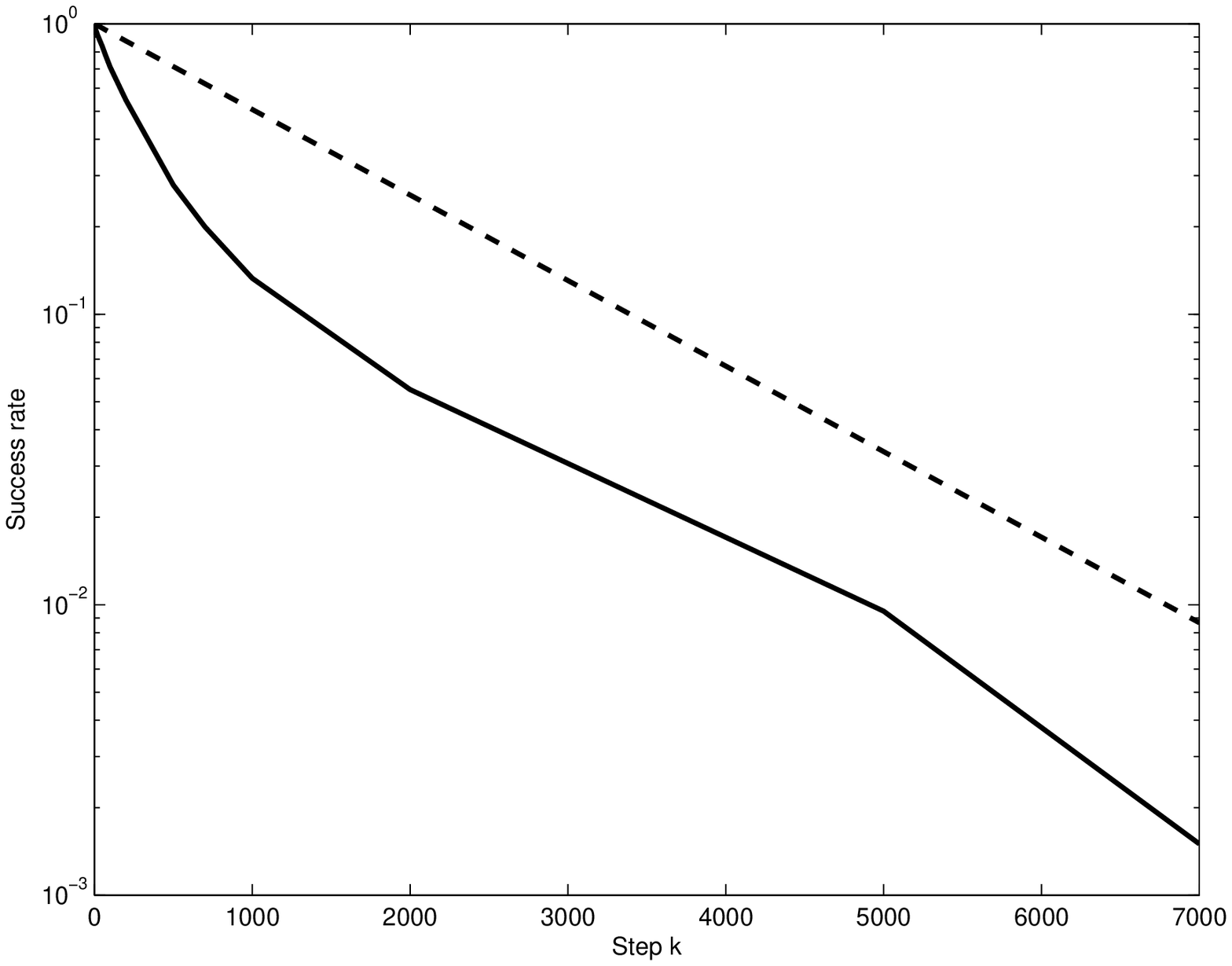}
}
\caption{Left panel: Success rate as a function of simulation step for
$J=1$ (dot-dash), $J=10$ (dotted), $J=100$ (solid) and $J=200$ (dashed).
Right panel: Logarithmic plot of
success rate for $J=100$ (solid) with bound of
Theorem \ref{thm:MT} (dashed).}
\label{fig:converge}
\end{figure}

To demonstrate the convergence of the algorithm, $2,000$ independent runs
of the algorithm, of $10,000$ steps each, were generated.
We then computed  the fraction of runs that found themselves in $\Theta^*(0.1)$
at different time points; for simplicity we refer to this fraction as the `success rate'.
The results for different values of $J$ are reported in the left panel of Fig.~\ref{fig:converge}.
It is clear that in all cases the success rate quickly settles to a steady
state value, suggesting that the algorithm has converged. Moreover, the steady
state success rate increases as $J$ increases. In the right panel of
Fig.~\ref{fig:converge} we concentrate on the case $J=100$ and plot
in a logarithmic scale the absolute value of the difference between the success
rate at different time points and the steady state success rate.
According to Theorem \ref{thm:MT}, one would expect this difference to decay
to $0$ geometrically at a rate $1-\frac{1}{M}$.
For comparison purposes, the corresponding curve for
the numerically estimated value $M=1475$, is also plotted on the figure.
The bound of Theorem \ref{thm:MT} indeed appears
to be valid, albeit, in this case, conservative.

\begin{figure}[h!]
\centerline{
  \includegraphics[width=0.5\hsize]{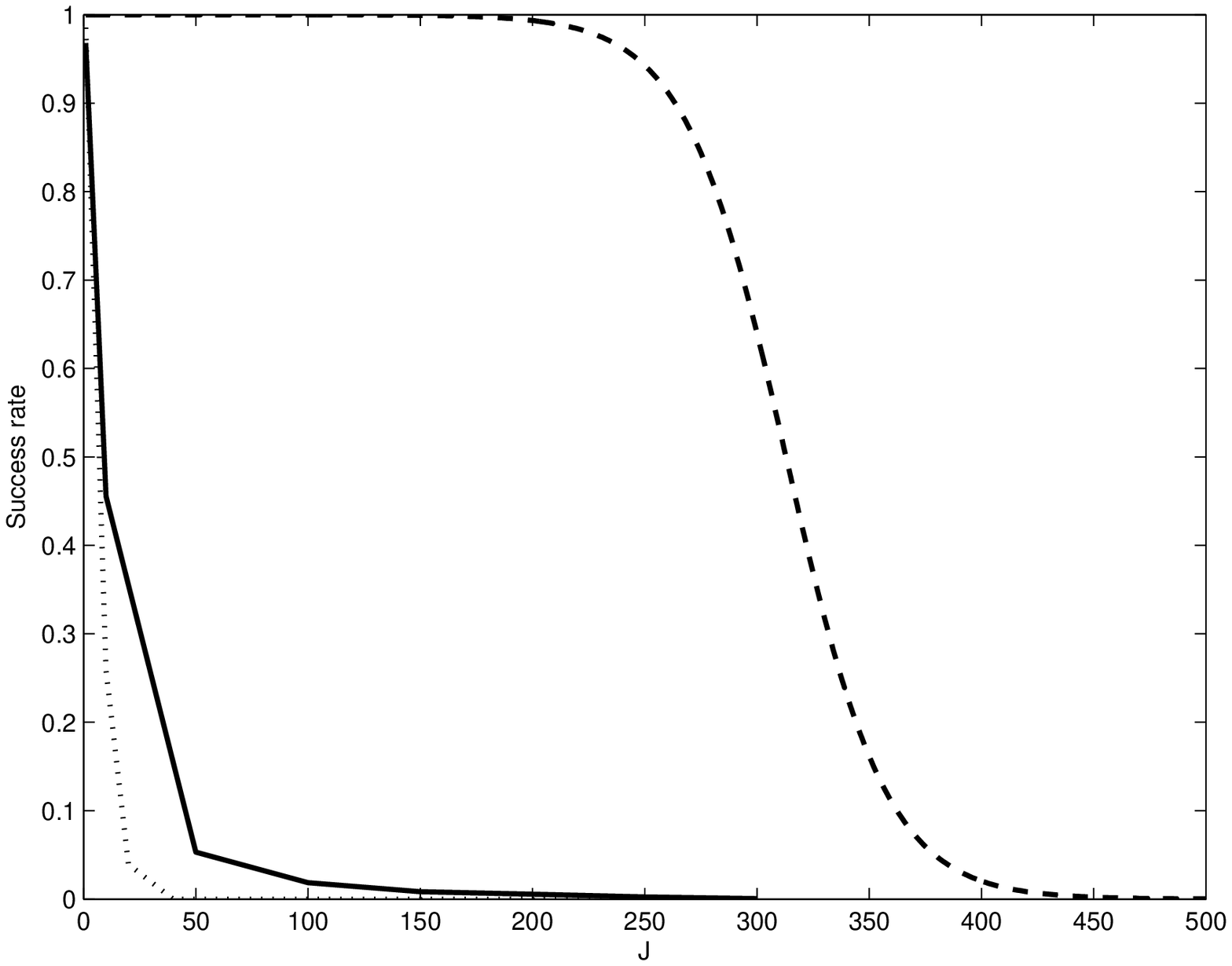}
  \includegraphics[width=0.5\hsize]{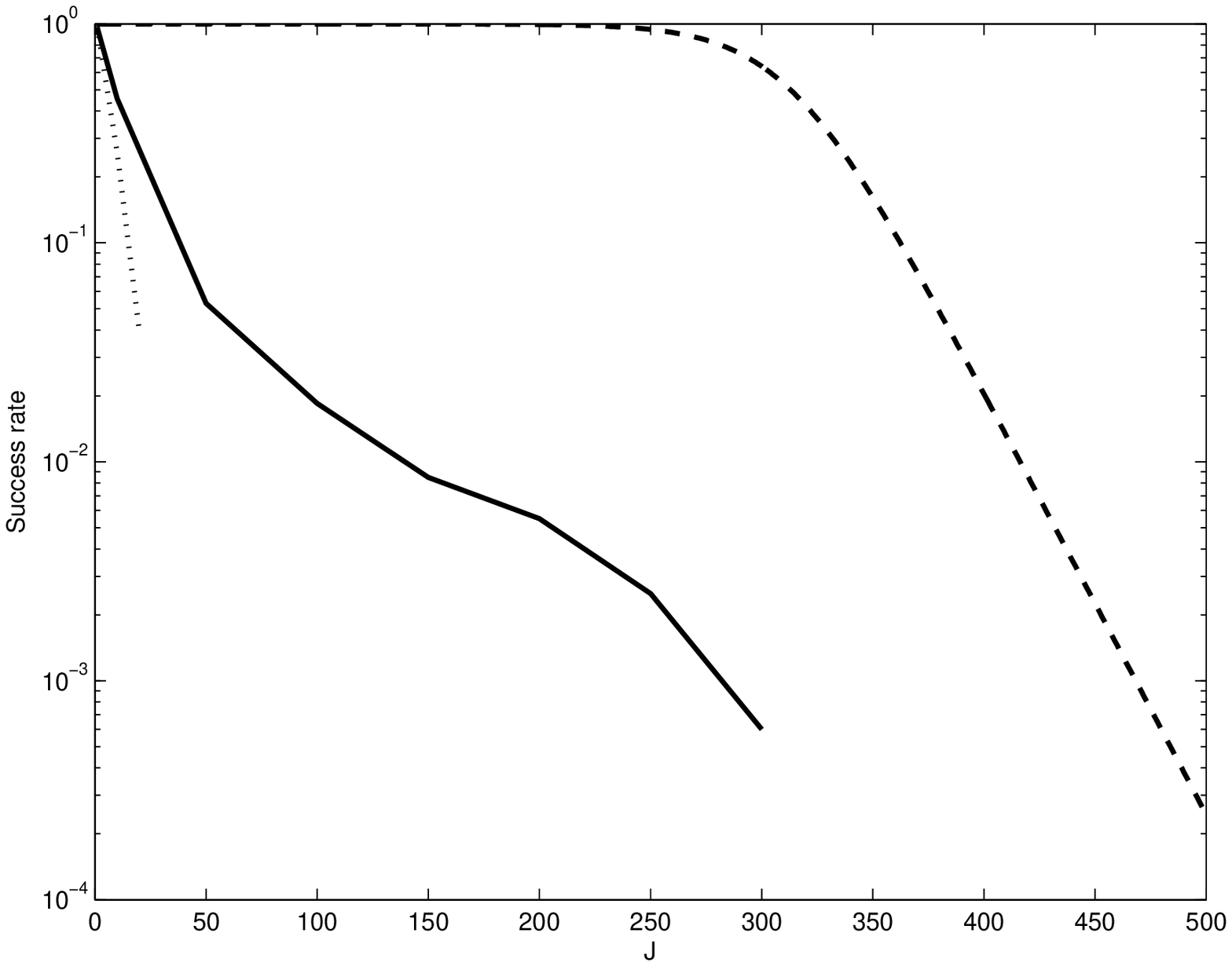}
}
\caption{Plot of the decay of $1$ minus the success rate as a function of the
exponent $J$ in linear (left panel) and logarithmic (right panel) scales.
Both plots show the empirical value based on the last state of $2,000$
independent runs of $10,000$ steps each (solid), the empirical
value based on the last $2,000$ states of a single $10,000$ step run
(dotted), and the theoretical bound (dashed).
$5,000$ independent runs were used for the case $J=300$, since the first
$2,000$ runs all ended \mbox{up in $\Theta^*(0.1)$ at step $10,000$.}\vspace{-0.5cm}}
\label{fig:sigmarate}
\end{figure}

\begin{figure}[t]
\centerline{
\includegraphics[width=0.5\hsize]{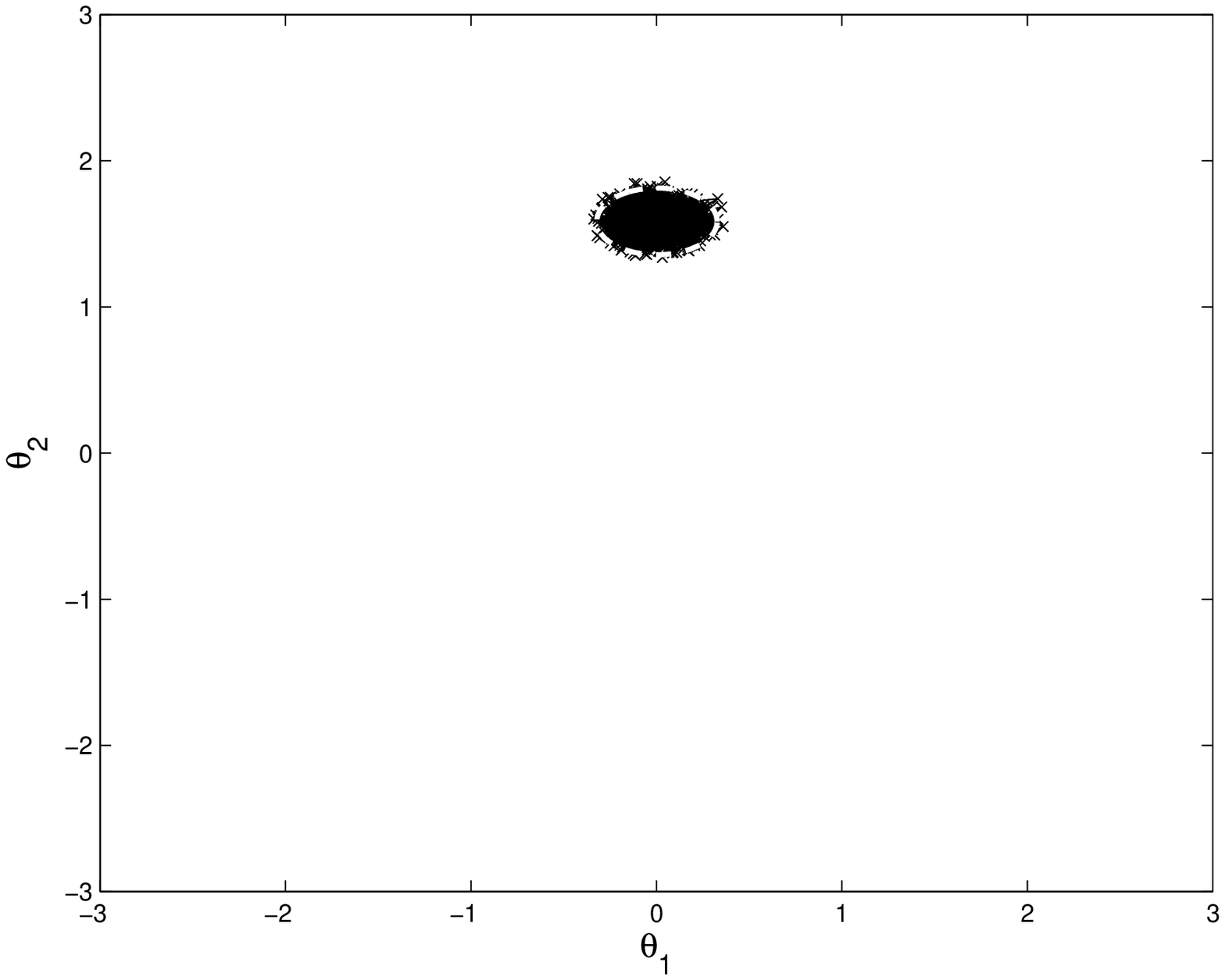}
\includegraphics[width=0.5\hsize]{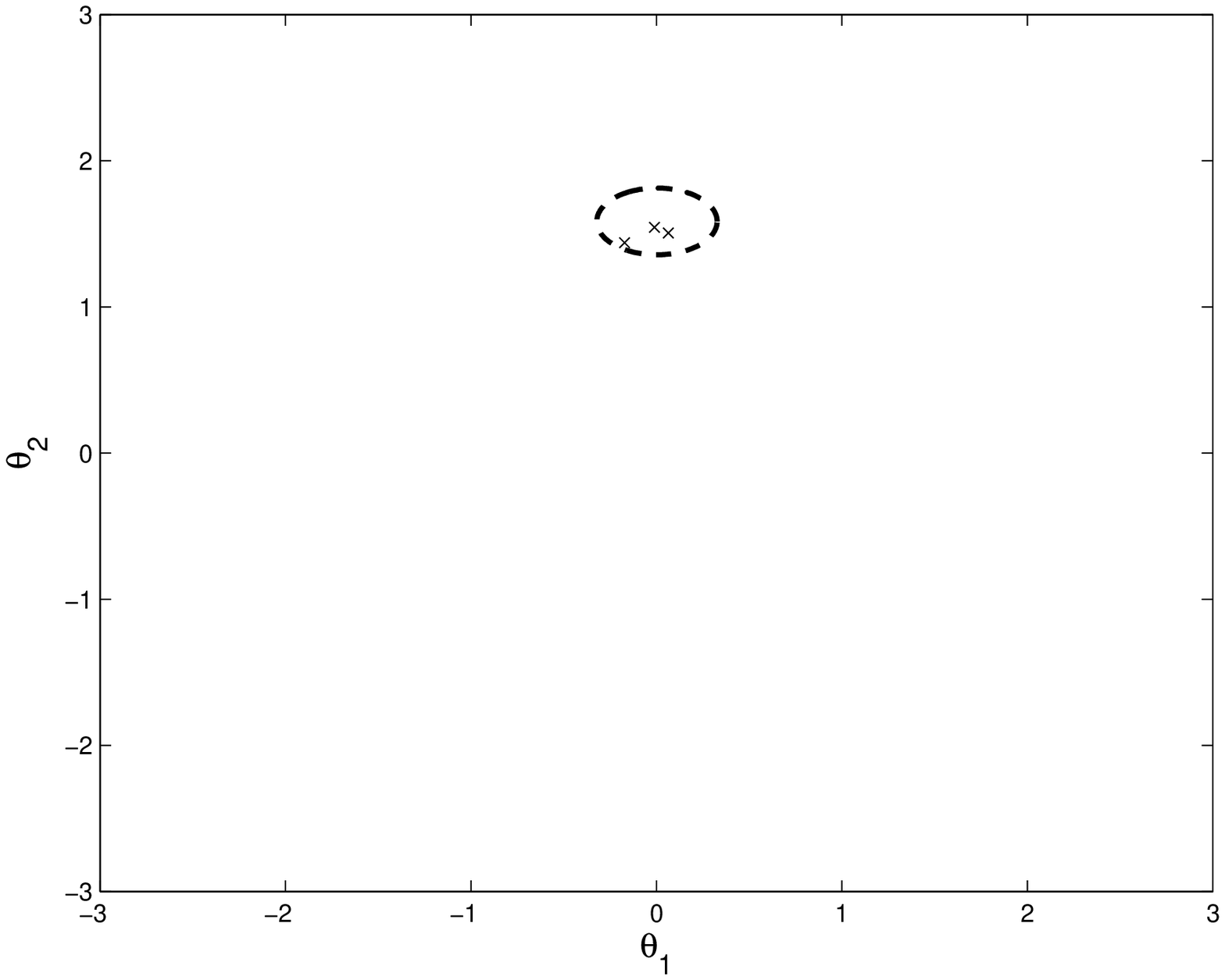}}
\centerline{
\includegraphics[width=0.5\hsize]{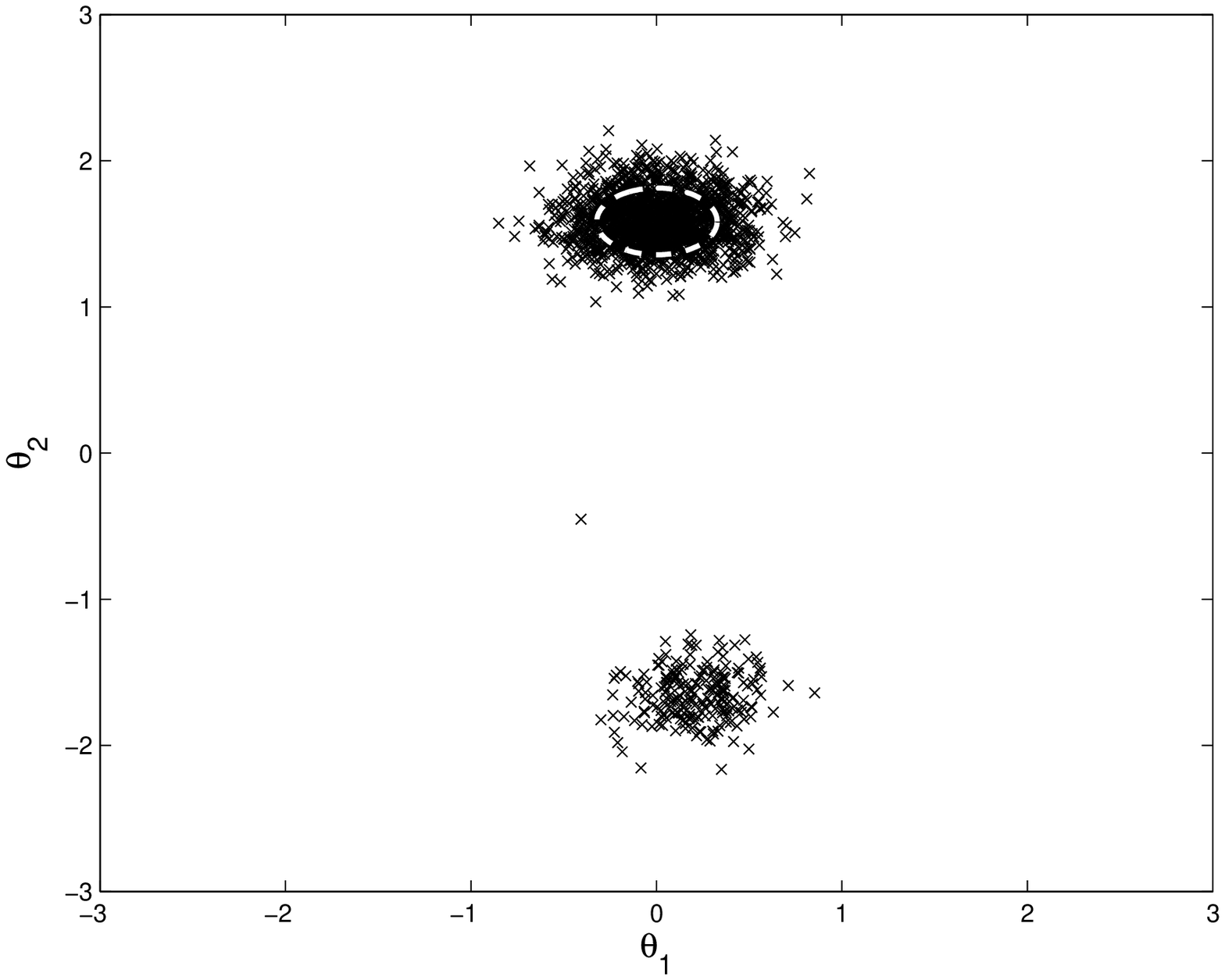}
\includegraphics[width=0.5\hsize]{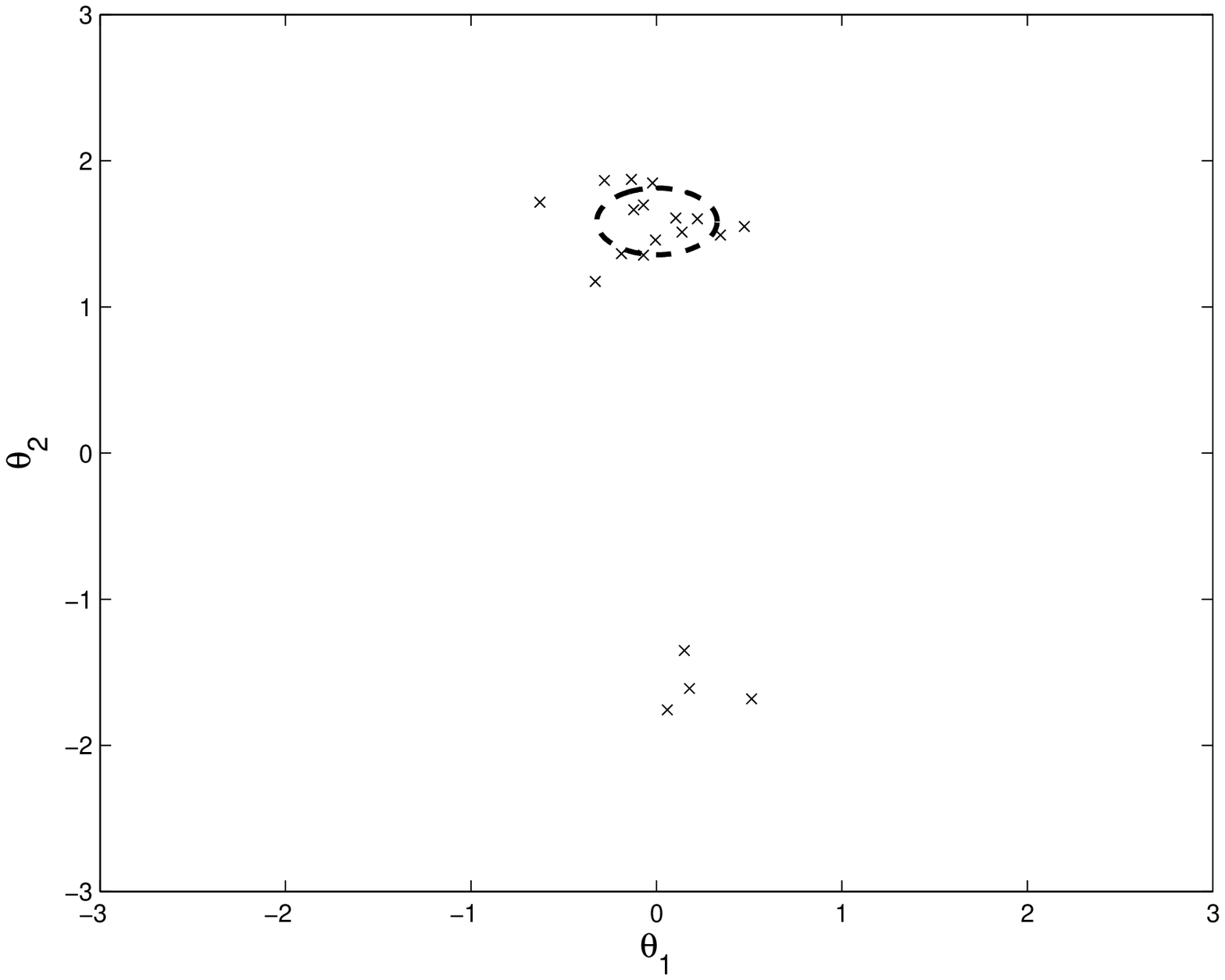}
}
\caption{Location of the the last state of $2,000$ independent runs of $10,000$
steps each (left column) and the last $2,000$ steps of a single $10,000$ step run
(right column) for $J=100$ (top row) and $J=10$ (bottom row). The set
$\Theta^*(0.1)$ is plotted as a dashed ellipse for comparison.}
\label{fig:tail}
\end{figure}

To demonstrate the bound of Proposition~\ref{prop:confidence_kR} the
steady state success rate as a function of the exponent $J$ is reported in Fig.~\ref{fig:sigmarate};
more precisely, the figure shows the decay of $1$ minus the steady state success rate
as a function of $J$ in linear and logarithmic scales.
The figure also shows the corresponding theoretical bound based on Proposition~\ref{prop:confidence_kR}.
Once again the bound appears to be valid.
Finally,  notice that although in this particular case the proposal
distribution  is an independent uniform distribution, the resulting states of
the chain are a sequence of dependent samples.
In Fig.~\ref{fig:sigmarate} we show the success rate estimated
using the last $2,000$ states of a single MCMC run of length $10,000$ (instead
of the last state of $2,000$ independent MCMC runs of length $10,000$ each).
It appears that the success rate increases much faster in this
case. This is due to the fact that the $2,000$ samples used are
now correlated. Figure~\ref{fig:tail} demonstrates this through a scatter
plots of the location of the $2,000$ states used to estimate the success rate
for $J=100$ in the two cases. While for both the $2,000$ independent runs
and the single run most of the states end up inside the set $\Theta^*(0.1)$ as
expected, it is apparent that the chain only moves three times in the last
$2,000$ steps of the single run; all other proposals are rejected. Plots
for the case $J=10$ are also included. Note that in this case points
near the second largest local maximum are also occasionally accepted.

\section{Comparison with other approaches to stochastic optimization}\label{sec:Disc}

In this section we attempt a comparison between the
computational features of the MCMC approach
with those of other state-of-the-art methods
for solving the stochastic programming problem (\ref{eq:exval})
with  finite-time  performance bounds.
Other methods are typically formulated under the assumption
that the distribution $P_{\bm{x}}(\cdot;\theta)$ does not depend on $\theta$.
In this case, $U$ becomes
\begin{equation}\label{eq:UPx}
U(\theta) = \int g(x,\theta)P_{\bm{x}}(dx).
\end{equation}
We stress from the beginning that a direct comparison of the computational
complexity of the different methods is not possible at this stage, since the
different methods rely on different assumptions, e.g.~some methods require
solving an additional optimization problem. Moreover, a satisfactory
complexity analysis for the MCMC approach is not yet available.
The comparison focuses on the number of samples of $\bm{x}$
required in each method  to obtain an approximate value optimizer
with imprecision $2\epsilon$, with confidence $\rho$,
in the optimization of (\ref{eq:UPx}).
In Table I we compare the growth rates of the total number of
samples of $\bm{x}$ required in each method to obtain the
desired optimization accuracy as a function of the parameters
of the problem.

In the approach of Shapiro \cite{Shapiro-Nemirovski-05,Shapiro-08},
$N$ independent samples $\bm{x}_1,\ldots,\bm{x}_N$, generated
according to $P_{\bm{x}}$, are used to construct the approximate
criterion
\begin{equation}\label{eq:uhat}
\hat{U}(\theta)=\frac{1}{N}\sum_{i=1}^N g(\bm{x}_i,\theta).
\end{equation}
It is shown in \cite{Shapiro-Nemirovski-05,Shapiro-08}
that if $\hat{\bm{\theta}}_S$ is an
approximate value optimizer of $\hat{U}$ with imprecision $\epsilon$
then $\hat{\bm{\theta}}_S$ is also an approximate value optimizer
of $U$ with imprecision $2\epsilon$ with probability at least $\rho$,
provided that $N$ is sufficiently high.
The growth rates of the required $N$,
reported in the first line of Table I, are based on
\cite[equation (3.9)]{Shapiro-08}. Notice that this is only a bound on the
samples required to construct $\hat{U}$.
It is argued in~\cite{Shapiro-08} that the optimization of $\hat U$
within  $\epsilon$ of optimality
can be carried out efficiently under convexity assumptions.
Nesterov \cite{Nesterov-05,Nesterov-Vial-08} presents a specific
approach for convex stochastic problems.
In this approach the samples
generated according to $P_{\bm{x}}$ are used to construct
an estimate of the optimizer of $U$ using
a stochastic sub-gradient algorithm.
The growth rates of the number of samples
required to obtain an approximate value optimizer of $U$
with imprecision $2\epsilon$ with probability at least $\rho$,
are reported in the second line of Table I and are based
on \cite[equation (14)]{Nesterov-Vial-08}.
Finally Vidyasagar~\cite{Vidyasagar-03,Vidyasagar-01} proposed a fully
randomized algorithm which, as mentioned earlier, is closely related to
the one presented in our work.
In Vidyasagar's approach one generates $N$ independent
samples $\bm{\theta}_1, \ldots, \bm{\theta}_N$ according to a `search'
distribution $P_{\bm{\theta}}$, which has support on $\Theta$, and
$M$ independent samples $\bm{x}_1, \ldots, \bm{x}_M$ according to $P_{\bm{x}}$,
and sets
\begin{equation}
\hat{\bm{\theta}}_V = \arg\min_{i=1, \ldots, N} \frac{1}{M}\sum_{j=1}^M g(\bm{x}_j,\bm{\theta}_i).
\label{eq:OptVid}
\end{equation}
Under minimal assumptions, close to our Assumption~\ref{assu:1}, it can be shown that if
\begin{equation}\label{eq:VidBound}
N \geq \frac{\mbox{$\displaystyle \log\frac{2}{1-\rho}$}}
{\mbox{$\displaystyle \log\frac{1}{1-\alpha}$}}
\quad \mbox{ and }\quad
M \geq \frac{1}{2\epsilon^2}\log\frac{4N}{1-\rho},
\end{equation}
then
\begin{equation}
P_{\bm{\theta}}(\{\theta \in \Theta \; | \;
U(\theta)>U(\hat{\bm{\theta}}_V)+\epsilon\}) \leq \alpha
\label{eq:measureVid}
\end{equation}
with probability at least $\rho$.
It is shown in \cite{Vidyasagar-01} that potentially tighter bounds
can be obtained if the family of functions $\{g(\cdot,\theta) \; | \; \theta \in \Theta\}$
has the UCEM property.
Notice that (\ref{eq:measureVid}) resembles (\ref{eq:DomOpt}) in our definition
of approximate domain optimizer.
The difference is that the measure of the set of points which are $\epsilon$ better
than the candidate optimizer is taken with respect to $P_{\bm{\theta}}$ in (\ref{eq:measureVid})
as opposed to the Lebesgue measure in (\ref{eq:DomOpt}).
If, and only if, $P_{\bm{\theta}}$ is chosen to be
the uniform  distribution over $\Theta$ then (\ref{eq:measureVid}) becomes
virtually equivalent to  (\ref{eq:DomOpt}).
In this case,  we can apply Theorem \ref{thm:1} and obtain the number of
samples required to obtain an approximate value optimizer.
By substituting $\alpha$ with the right-hand side of (\ref{eq:alphaboundR}) in
(\ref{eq:VidBound}) 
we obtain that if
\begin{equation}\label{eq:VidBoundValue}
N \geq \left(\frac{LR}{\epsilon}\right)^n\log\frac{2}{1-\rho}
\quad\mbox{ and }\quad
M \geq \frac{1}{2\epsilon^2} \left[\log\frac{4}{1-\rho} + \log\log\frac{2}{1-\rho}+n\log\frac{LR}{\epsilon}\right]
\end{equation}
then $\hat\theta_V$ is an approximate value optimizer of $U$
with imprecision $2\epsilon$ with probability at least $\rho$.
Notice that now the number of samples on $\Theta$
turns out to be exponential in $n$.

\begin{table}\label{tab:growth}
\centering\scriptsize
\setlength\minrowclearance{2pt}
\hspace*{-0.2cm}\begin{tabular}{lcccccccccc}
\toprule
                                                                                    &\mbox{ } &$\displaystyle \epsilon$                                   &\mbox{ } & $\displaystyle \rho\mbox{ or }\sigma$           &\mbox{ }& $\displaystyle LR$             &\mbox{ } &$\displaystyle n$ &\mbox{ } & problem\\\otoprule
Shapiro  \cite{Shapiro-Nemirovski-05,Shapiro-08}                                    &\mbox{ } &$\displaystyle \frac{1}{\epsilon^2}\log\frac{1}{\epsilon}$ &\mbox{ } & $\displaystyle \log\frac{1}{1-\rho}$    &\mbox{ }& $\displaystyle (LR)^2\log LR$  &\mbox{ } &$\displaystyle n$ &\mbox{ } & convex\\\midrule
Nesterov \cite{Nesterov-05,Nesterov-Vial-08}                                        &\mbox{ } &$\displaystyle \frac{1}{\epsilon^4}$                       &\mbox{ } & $\displaystyle \log\frac{1}{1-\rho}$    &\mbox{ }& $\displaystyle (LR)^2$         &\mbox{ } &$\displaystyle -$ &\mbox{ } & convex\\\midrule
Vidyasagar (\ref{eq:VidBoundValue})                                                 &\mbox{ } &$\displaystyle \frac{1}{\epsilon^2}\log\frac{1}{\epsilon}$ &\mbox{ } & $\displaystyle \log\frac{1}{1-\rho}$    &\mbox{ }& $\displaystyle \log LR$        &\mbox{ } &$\displaystyle n$ &\mbox{ } & general\\\midrule
\rowcolor[gray]{.8} MCMC (\ref{eq:JboundvalueR}) {\sc [per iteration]}                &\mbox{ } & $\displaystyle \frac{1}{\epsilon}\log\frac{1}{\epsilon}$  &\mbox{ } & $\displaystyle \log\frac{1}{1-\sigma}$  &\mbox{ }& $\displaystyle \log LR$        &\mbox{ } &$\displaystyle n$ &\mbox{ } & general\\
\bottomrule
\end{tabular}
\begin{minipage}{\columnwidth}
\mbox{ }\\ \mbox{ }\\
{\sc TABLE I:} Growth rates of the number of samples of $\bm{x}$ required to obtain an approximate value optimizer
with imprecision $2\epsilon$
of $U$, given by (\ref{eq:UPx}), with probability $\rho$.
In the case of MCMC, the entries of the table represent the number of samples of $\bm{x}$
required to perform one iteration of the algorithm.
\end{minipage}
\end{table}

In the last row of the table, we have included
the growth rates of (\ref{eq:JboundvalueR}), which is
the number of samples of $\bm{x}$ which must be generated
at each iteration of Algorithm II (which coincides with the exponent $J$).
In this case, the total number of required samples of $\bm{x}$ is $J$
times the number of iterations required to achieve the desired reduction of
$\|P_{\bm{\theta}_k}  - \pi(\,\cdot\,\,;J,\delta)\|_{\mbox{\tiny TV}}$.
Hence, the entries of the last row represent
a lower bound, or the `base-line' growth rates,
of the total number of required samples in the MCMC approach.
In this case, the confidence is
$\rho = \sigma - \|P_{\bm{\theta}_k}  - \pi(\,\cdot\,\,;J,\delta)\|_{\mbox{\tiny TV}}$.
Hence, since $\sigma>\rho$, it is sensible to consider the growth rate with respect to  $\sigma$ instead of $\rho$.
By comparing the different entries in the table, we notice that (\ref{eq:JboundvalueR})
grows slower than or at the same rate as the other bounds.
Overall, the comparison reveals that in principle there is scope
for obtaining  MCMC algorithm which, in terms of numbers of required samples of $\bm{x}$,
have a computational cost comparable to those of the other approaches.
Here we present a preliminary complexity analysis
which shows that the introduction of a cooling schedule
would eventually lead to efficient algorithms.
Notice that in Section \ref{sec:PCA} we
considered a constant schedule ($J_k=J$).
Here, we assume that $J_k$ takes integer values starting with
$J_1=1$ and ending with $J_k=J$, where $J$
is the smallest integer which satisfies either
(\ref{eq:Jbound}) or (\ref{eq:JboundvalueR}).
In the earliest works on simulated annealing,
the logarithmic schedule of the type
$J_k = \lfloor \log k \rfloor + 1$
was often adopted \cite{Haario-Saksman-91,Gelfand-Mitter-91,Andrieu-et-al-01}.
Here, we are interested in counting the total number of iterations required
to complete the cooling schedule when $J$ is given by (\ref{eq:Jbound}) or (\ref{eq:JboundvalueR}).
Let $K_i$ denote the number of
iterations in which $J_k=i$ for each $i=1,2,\dots,J$.
Hence, the total number of iterations
is $\sum_{i=1}^J K_i$  and, in Algorithm II,
the total number of required samples of $\bm{x}$ would be $\sum_{i=1}^J iK_i$.
For the logarithmic schedule $J_k = \lfloor \log k \rfloor + 1$
we have $K_i = \lfloor e^{i}\rfloor- \lfloor e^{i-1}\rfloor$.
Hence we obtain
$$
\sum_{i=1}^J K_i =
\sum_{i=1}^J \lfloor e^{i}\rfloor- \lfloor e^{i-1}\rfloor
\approx  e^{J}-1\, ,\qquad
\sum_{i=1}^J iK_i  \approx Je^{J} - \frac{1-e^J}{1-e}.
$$
In this case the number of iterations  turns out to be exponential in $\frac{1}{\epsilon}$.
Hence, a logarithmic schedule is not sufficient to obtain efficient algorithms.

In more recent works  \cite{Tsallis-Stariolo-96,Locatelli-00,Rubenthaler-et-al-09},
the faster algebraic schedule of the type
$J_k = \lfloor k^a\rfloor$, with $a>0$, has been considered.
It is shown in \cite{Tsallis-Stariolo-96,Locatelli-00,Rubenthaler-et-al-09}
that the choice of the faster algebraic schedule
requires a sophisticated design of the proposal distribution.
For the algebraic schedule  $J_k = \lfloor k^a\rfloor$
we have $K_i = \lfloor (i+1)^\frac{1}{a}\rfloor - \lfloor i^\frac{1}{a}\rfloor$.
Hence we obtain
\begin{eqnarray*}
\sum_{i=1}^J K_i &=&
\sum_{i=1}^J \lfloor (i+1)^\frac{1}{a}\rfloor - \lfloor i^\frac{1}{a}\rfloor
\approx  \left(J+1\right)^{\frac{1}{a}} - 1\, ,\\
\sum_{i=1}^J iK_i
&\approx&
 J\left(J+1\right)^{\frac{1}{a}} - \sum_{i=0}^{J-1}\left(i +1 \right)^\frac{1}{a}.
\end{eqnarray*}
In this case the number of iterations grows as $J^\frac{1}{a}$.
Hence, in the case of approximate domain optimization,
where $J$ is given by (\ref{eq:Jbound}),
the number of iterations grows as
$(\frac{1}{\epsilon})^{\frac{1}{a}}$, $(\log\frac{1}{\alpha})^{\frac{1}{a}}$
and $(\log\frac{1}{1-\sigma})^{\frac{1}{a}}$.
In the case of approximate value optimization,
where $J$ is given by (\ref{eq:JboundvalueR}),
the number of iterations grows as
$(\frac{1}{\epsilon}\log\frac{1}{\epsilon})^{\frac{1}{a}}$,
$(\log\frac{1}{1-\sigma})^{\frac{1}{a}}$,
$(\log LR)^{\frac{1}{a}}$ and  $n^{\frac{1}{a}}$.
In Algorithm II, the total number of samples of $\bm{x}$ is given by
the number of iterations multiplied by $J$.
Hence, in the case of approximate domain optimization, it
grows as
$(\frac{1}{\epsilon})^{1+\frac{1}{a}}$, $(\log\frac{1}{\alpha})^{1+\frac{1}{a}}$
and $(\log\frac{1}{1-\sigma})^{1+\frac{1}{a}}$.
In the case of approximate value optimization, it grows as
$(\frac{1}{\epsilon}\log\frac{1}{\epsilon})^{1+\frac{1}{a}}$,
$(\log\frac{1}{1-\sigma})^{1+\frac{1}{a}}$,
$(\log LR)^{1+\frac{1}{a}}$ and $n^{1+\frac{1}{a}}$
(notice that, as $a$ increases, the growth rates
approach the entries of the last row of Table I).
Hence, an algebraic schedule leads to algorithms with polynomial growth rates.
The convergence analysis of Algorithms I and II with an additional algebraic
cooling schedule goes beyond the scope of this paper.
Here, we limit ourselves to pointing out that
the choice of a target distribution $\pi(\,\cdot\,\,;J,\delta)$
with a finite $J$ implies that the cooling
schedule $\{J_k\}_{k=1,2,\dots}$ can be chosen to be a sequence that
takes only a finite set of values.
In turn, this fact should make the study  of convergence
of $P_{\bm{\theta}_k}$ to $\pi(\,\cdot\,\,;J,\delta)$ in total variation
distance easier than the study of asymptotic convergence of
$P_{\bm{\theta}_k}$  to the zero-temperature distribution $\pi_{\infty}$.


\section{Conclusions}\label{sec:end}

In this paper, we have introduced a novel approach
for obtaining  rigorous finite-time guarantees
on the performance of MCMC algorithms
in the optimization of
functions of continuous variables.
In particular we have established the values of the
the temperature parameter in the target distribution
which allow one to reach a solution,
which is within the desired level of approximation
with the desired confidence,  in a finite number of steps.
Our work was motivated by the MCMC algorithm (Algorithm II),
introduced in~\cite{Muller-99,Doucet-et-al-02,Muller-et-al-04},
for solving stochastic optimization problems.
On the basis of our results, we were able to
obtain the `base-line' computational complexity
of the MCMC approach and to  perform an initial
assessment of the computational complexity
of MCMC algorithms. It has been shown that
MCMC algorithms with an algebraic cooling schedule would
have polynomial complexity bounds comparable with those
of other state-of-the-art
methods for solving stochastic optimization problems.
Conditions for asymptotic convergence of
simulated annealing algorithms with an algebraic cooling schedule
have already been reported in the literature \cite{Tsallis-Stariolo-96,Locatelli-00,Rubenthaler-et-al-09}.
Our results enable novel research on the development of
efficient MCMC algorithms for the solution of stochastic programming
problems with rigorous finite-time guarantees.
Finally, we would like to point out that the results
presented in this work do not apply to the  MCMC approach only
but do apply also to other sampling methods which can implement the idea of
simulated annealing \cite{DelMoral-Miclo-03,DelMoral-et-al-06}.

\begin{center}
{\sc Acknowledgments}
\end{center}
Work supported by EPSRC, Grant EP/C014006/1, and by the European Commission under
projects HYGEIA FP6-NEST-4995 and iFly FP6-TREN-037180.

\appendix
In order to prove Theorem~\ref{thm:1} we first need to prove a preliminary technical result.
\begin{lemma}\label{prop:lemma}
For $A \subseteq \Re^n$ and $\theta \in \Re^n$ let
$
d(\theta,A):=\inf_{\theta'\in A}\|\theta-\theta'\|.$
Then, for any $\beta \geq 0$
$$
d^*_\beta:=\sup_{\tiny\begin{array}{c} A\subseteq \Re^n
\\
\lambda(A)\leq \beta\end{array}}
\sup_{\theta \in A}d(\theta,A^c)=
\frac{1}{\sqrt{\pi}}\left[\frac{n}{2}\Gamma\left(\frac{n}{2}\right)
\right]^{\frac{1}{n}}\beta^{\frac{1}{n}}
$$
where $A^c$ denotes the complement of $A$ in $\Re^n$.
\end{lemma}
In the above Lemma the inner supremum determines the
points in the set $A$ whose distance from the complement of $A$ is the
largest; loosely speaking the points that lie the furthest from the
boundary of $A$ or the deepest in the interior of $A$. The outer
supremum then maximizes this distance over all sets $A$ whose Lebesgue
measure is bounded by $\beta$.

\begin{proof_text}
We show that the optimizers for the outer supremum are 2-norm
balls in $\Re^n$; then, the inner supremum is achieved
by the center of the ball.

Consider any set
$A \subseteq \Re^n$ with $\lambda(A)\leq \beta$ and let
\begin{equation}
d_A:=\sup_{\theta \in A}d(\theta,A^c). \label{eq:sup1}
\end{equation}
Since $\lambda(A)<\infty$ , the supremum~\eqref{eq:sup1}
is achieved by some $\theta_A \in A$. Without loss of generality we
can assume that the set $A$ is closed; if not, taking its closure will
not affect its Lebesgue measure and will lead to the same value for $d_A$.
Let $B(\theta,d)$ denote the 2-norm ball
with center in $\theta$ and radius $d$.
Notice that by construction $B(\theta_A,d_A) \subseteq A$ and
therefore $\lambda(B(\theta_A,d_A)) \leq \lambda(A)\leq \beta$.
Moreover,
$$
\sup_{\theta \in B(\theta_A,d_A)}d(\theta,B(\theta_A,d_A)^c)=d_A
$$
achieved at the center, $\theta_A$, of the ball.
In summary, for any $A \subseteq \Re^n$ with $\lambda(A)\leq\beta$
one can find a 2-norm ball of measure at most $\beta$ which achieves
$\sup_{\theta \in A}d(\theta,A^c)$.

Therefore,
\begin{align*}
d^*_\beta :&= \sup_{\tiny\begin{array}{c} A\subseteq \Re^n
\\
\lambda(A)\leq \beta\end{array}}
\sup_{\theta \in A}d(\theta,A^c)\\
& =\sup_{\tiny\begin{array}{c}(\theta',r)\in \Re^n\times \Re_+\\
\lambda(B(\theta',r))\leq \beta\end{array}} \sup_{\theta \in
  B(\theta',r)}d(\theta,B(\theta',r)^c)\\
& = \sup_{\tiny\begin{array}{c}r \geq 0\\
\lambda(B(0,r))\leq \beta\end{array}} \sup_{\theta \in
  B(0,r)}d(\theta,B(0,r)^c)\\
& = \sup_{\tiny\begin{array}{c}r \geq 0\\
\lambda(B(0,r))\leq \beta\end{array}}r\\
& =
\frac{1}{\sqrt{\pi}}\left[\frac{n}{2}\Gamma\left(\frac{n}{2}\right)
\right]^{\frac{1}{n}}\beta^{\frac{1}{n}}.
\end{align*}
In the above derivation, the last equality is obtained by recalling that
$\displaystyle \lambda(B(\theta,r)) = \frac{2\pi^{\frac{n}{2}}}{n\Gamma(\frac{n}{2})}r^n$.
\end{proof_text}
We are now in a position to prove Theorem~\ref{thm:1}.

\begin{proof_ref}{Theorem~\ref{thm:1}}
Let $\hat{\Theta}(\theta,\epsilon):=\{\theta' \in \Theta \; | \;U(\theta')>U(\theta)+\epsilon\}$
and recall that by definition
$$
\lambda(\hat{\Theta}(\theta,\epsilon)) \leq \alpha\, \lambda(\Theta)\, .
$$
Take any $\tilde{\theta} \in \Theta$.
Then either
$U(\tilde{\theta})\leq U(\theta)+\epsilon$  or $\tilde{\theta} \in \hat{\Theta}(\theta,\epsilon)$.
In the former case there is nothing to prove.
In the latter case, according to Lemma \ref{prop:lemma},
we have:
$$
d(\tilde\theta,\hat{\Theta}(\theta,\epsilon)^c)
\leq
\frac{1}{\sqrt{\pi}}\left[\frac{n}{2}\Gamma\left(\frac{n}{2}\right)
\right]^{\frac{1}{n}}[\alpha
  \lambda(\Theta)]^{\frac{1}{n}}.
$$
Since $\Theta$ is compact and $U$ is continuous, the set
$\hat{\Theta}(\theta,\epsilon)^c$ is closed and therefore
there exists  $\bar\theta \in \hat{\Theta}(\theta,\epsilon)^c$
such that
\[
\| \tilde\theta-\bar\theta\| \leq
  \frac{1}{\sqrt{\pi}}\left[\frac{n}{2}\Gamma\left(\frac{n}{2}\right)
  \right]^{\frac{1}{n}}[\alpha \lambda(\Theta)]^{\frac{1}{n}}.
\]
Moreover, since $U$ is Lipschitz,
$
|U(\tilde\theta)-U(\bar\theta)| \leq L \| \tilde\theta-\bar\theta\|$.
Since $\bar\theta \in \hat{\Theta}(\hat{\theta},\epsilon)^c$, we have that
$U(\bar\theta)\leq U(\theta)+\epsilon$, and therefore
\begin{equation}\label{eq:sup tilde}
U(\tilde\theta) \leq U(\theta) + \epsilon +
\frac{L}{\sqrt{\pi}}
\left[\frac{n}{2}\Gamma\left(\frac{n}{2}\right)\right]^{\frac{1}{n}}
[\alpha \lambda(\Theta)]^{\frac{1}{n}}.
\end{equation}
The claim follows since $\tilde\theta$
is arbitrary in $\Theta$ and satisfies either
$U(\tilde{\theta})\leq U(\theta)+\epsilon$ or (\ref{eq:sup tilde}).
\end{proof_ref}

\begin{proof_ref}{Theorem~\ref{thm:bound}}
Let ${\bar\alpha}\in (0, 1]$ and $\rho \in (0,1]$ be given numbers.
To simplify the notation, let
$U_\delta(\theta):=U(\theta)+\delta$ and
let $\pi_\delta$ be a normalized measure such that
$\pi_\delta(d\theta) \propto U_\delta(\theta)\lambda(d\theta)$,
i.e.~$\pi_\delta(d\theta):=\pi(d\theta;1,\delta)$.
In the first part of the
proof we establish a lower bound on
$$
\pi\left(\{\theta\in\Theta \; | \; \pi_\delta(\{\theta'\in\Theta\; | \;\rho\,U_\delta(\theta') > U_\delta(\theta)\}) \leq {\bar\alpha}\};J,\delta\right)\, .
$$

Let
$y_{\bar\alpha} := \inf \{y\; | \; \pi_\delta (\{\theta\in\Theta\; | \; U_\delta(\theta)\leq y\}) \geq  1-{\bar\alpha} \}$.
To start with we  show that the set
$\{\theta\in\Theta\; | \; \pi_\delta(\{\theta'\in\Theta\; | \; \rho\, U_\delta(\theta')> U_\delta(\theta) \})\leq{\bar\alpha} \}$
coincides with $\{ \theta\in\Theta\; | \;  U_\delta(\theta) \geq \rho\, y_{\bar\alpha} \}$.
Notice that the quantity $\pi_\delta (\{\theta\in\Theta\; | \; U_\delta(\theta)\leq y\})$
is a non decreasing right continuous function of $y$ because it has the form of a distribution function
(see e.g.~\cite[p.~162]{Gnedenko-68}, see also \cite[Lemma 11.1]{Vidyasagar-03}).
Therefore we have
$\pi_\delta (\{\theta\in\Theta\; | \; U_\delta(\theta)\leq y_{\bar\alpha}\}) \geq 1-{\bar\alpha}$ and
$$
y\geq \rho\, y_{\bar\alpha}\,\Rightarrow\,
\pi_\delta (\{\theta'\in\Theta\; | \; \rho\,  U_\delta(\theta')\leq  y\}) \geq 1-{\bar\alpha}\,
\Rightarrow\, \pi_\delta  (\{\theta'\in\Theta\; | \; \rho\, U_\delta(\theta') >  y\})\leq {\bar\alpha}\, .
$$
Moreover,
$$
y < \rho\, y_{\bar\alpha}\, \Rightarrow\,
\pi_\delta (\{\theta'\in\Theta\; | \; \rho\,  U_\delta(\theta')\leq y\}) < 1-{\bar\alpha}\, \Rightarrow\,
\pi_\delta (\{\theta'\in\Theta\; | \; \rho\,  U_\delta(\theta') > y\}) > {\bar\alpha}\,
$$
and taking the contrapositive one obtains
$$
\pi_\delta  (\{\theta'\in\Theta\; | \; \rho\, U_\delta(\theta') > y\})\leq{\bar\alpha}\quad\Rightarrow\quad
y\geq \rho\, y_{\bar\alpha}\, .
$$
Therefore $\{ \theta\in\Theta\; | \; U_\delta(\theta) \geq\rho\, y_{\bar\alpha} \} =
\{\theta\in\Theta\; | \; \pi_\delta(\{\theta'\in\Theta\; | \;\rho\,  U_\delta(\theta')> U_\delta(\theta) \})\leq {\bar\alpha} \}$.

We now derive a lower bound on
$\pi\left(\{ \theta\in\Theta\; | \;  U_\delta(\theta)\geq \rho\, y_{\bar\alpha} \};J,\delta\right)$.
Let us introduce the notation
$A_{\bar\alpha}:=\{\theta\in\Theta\; | \; U_\delta(\theta)<y_{\bar\alpha}\}$,
$\bar A_{\bar\alpha}:=\{\theta\in\Theta\; | \;  U_\delta(\theta)\geq y_{\bar\alpha}\}$,
$B_{{\bar\alpha},\rho}:=\{\theta\in\Theta\; | \; U_\delta(\theta)<\rho\, y_{\bar\alpha}\}$
and
$\bar B_{{\bar\alpha},\rho}:=\{\theta\in\Theta\; | \; U_\delta(\theta)\geq \rho\, y_{\bar\alpha}\}$.
Notice that
$B_{{\bar\alpha},\rho}\subseteq A_{\bar\alpha}$ and
$\bar A_{\bar\alpha} \subseteq \bar B_{{\bar\alpha},\rho}$.
The quantity
$\pi_\delta (\{\theta\in\Theta\; | \;  U_\delta(\theta)< y\})$
as a function of $y$ is the left continuous version of
$\pi_\delta (\{\theta\in\Theta\; | \; U_\delta(\theta)\leq y\})$~\cite[p.~162]{Gnedenko-68}.
Hence, the definition of $y_{\bar\alpha}$
implies $\pi_\delta(A_{\bar\alpha})\leq 1-{\bar\alpha}$ and $\pi_\delta(\bar A_{\bar\alpha})\geq{\bar\alpha}$.
Notice that
\begin{align*}
\pi_\delta(A_{\bar\alpha})\leq 1-{\bar\alpha} \quad&\Rightarrow\quad
\frac{\delta\lambda(A_{\bar\alpha})} {\left[\int_{\Theta} U_\delta(\theta) \lambda(d\theta)\right]}
\leq 1-{\bar\alpha}\qquad\mbox{ because } U(\theta)\geq0\; \forall \theta\, ,\\
\pi_\delta(\bar A_{\bar\alpha}) \geq {\bar\alpha} \quad
&\Rightarrow\quad
\frac{(1+\delta)\lambda(\bar A_{\bar\alpha})}{\left[\int_{\Theta} U_\delta(\theta)\lambda(d\theta)\right]}\geq{\bar\alpha}
\qquad\qquad\mbox{because } U(\theta)\leq 1\; \forall \theta\, .
\end{align*}
Hence, ${\lambda(\bar A_{\bar\alpha})}>0$ and
$$
\frac{\lambda(A_{\bar\alpha})}{\lambda(\bar A_{\bar\alpha})}\leq
\frac{1-{\bar\alpha}}{{\bar\alpha}}\frac{1+\delta}{\delta}\, .
$$
Notice that ${\lambda(\bar A_{\bar\alpha})}>0$ implies ${\lambda(\bar B_{{\bar\alpha},\rho})}>0$.
We obtain
\begin{eqnarray*}
\pi\left(
\{\theta\in\Theta\; | \;  U_\delta(\theta) \geq \rho\,y_{\bar\alpha}\};J,\delta\right)
&=&
\pi\left(
\bar B_{{\bar\alpha},\rho};J,\delta\right)
=
\frac
{
\int_{\bar B_{{\bar\alpha},\rho}}
  U_\delta(\theta)^J\lambda(d\theta)}
{
\int_{\Theta}
U_\delta(\theta)^J\lambda(d\theta)}\\
&=&
\frac
{
\int_{\bar B_{{\bar\alpha},\rho}}
  U_\delta(\theta)^J\lambda(d\theta)}
{
\int_{B_{{\bar\alpha},\rho}}
U_\delta(\theta)^J\lambda(d\theta)+
\int_{\bar B_{{\bar\alpha},\rho}}
U_\delta(\theta)^J\lambda(d\theta)}\\
&=&
\frac{1}
{\mbox{$\displaystyle 1 +
\frac{\int_{ B_{{\bar\alpha},\rho}}
  U_\delta(\theta)^J\lambda(d\theta)}
{\int_{\bar B_{{\bar\alpha},\rho}}
  U_\delta(\theta)^J\lambda(d\theta)} $}}\\
&\geq&
\frac{1}
{\mbox{$\displaystyle 1 +
\frac{\int_{B_{{\bar\alpha},\rho}}
  U_\delta(\theta)^J\lambda(d\theta)}{\int_{\bar A_{\bar\alpha}}
  U_\delta(\theta)^J\lambda(d\theta)} $}}\\
&\geq&
\frac{1}
{\mbox{$\displaystyle 1 +
 \frac{\rho^{\,
     J}y_{\bar\alpha}^J}{y_{\bar\alpha}^J}
 \frac{\lambda(B_{{\bar\alpha},\rho})}{\lambda(\bar A_{\bar\alpha})}$}}\\
&\geq& \frac{1}
{\mbox{$\displaystyle 1 +
\rho^{\,J}\frac{\lambda(A_{\bar\alpha})}{\lambda(\bar A_{\bar\alpha})}$}}\\
&\geq&
\frac{1}
{\mbox{$\displaystyle 1 +
\rho^{\,J}\frac{1-{\bar\alpha}}{{\bar\alpha}}\frac{1+\delta}{\delta}$}}\, .
\end{eqnarray*}
Since
$\{ \theta\in\Theta\; | \;  U_\delta(\theta) \geq \rho\, y_{\bar\alpha} \}=
\{\theta\in\Theta\; | \; \pi_\delta(\{\theta'\in\Theta\; | \; \rho\, U_\delta(\theta')> U_\delta(\theta) \})\leq{\bar\alpha} \}$
the first part of the proof is complete.

In the second part of the proof we show that the set
$\{\theta\in\Theta\; | \; \pi_\delta(\{\theta'\in\Theta\; | \; \rho\, U_\delta(\theta')> U_\delta(\theta) \})\leq{\bar\alpha} \}$
is contained in the set of approximate domain
optimizers of $U$ with value imprecision
$\tilde\epsilon := (\rho^{-1}-1)(1+\delta)$ and residual domain
$\tilde\alpha := \frac{1+\delta}{\tilde\epsilon + \delta}\,{\bar\alpha}$.
Hence, we show that
\begin{eqnarray*}
\lefteqn{\{\theta\in\Theta\; | \; \pi_\delta(\{\theta'\in\Theta\; | \; \rho\, U_\delta(\theta')> U_\delta(\theta) \})\leq{\bar\alpha} \}
\subseteq}\\
& &\{\theta\in\Theta\; | \; \lambda(\{\theta'\in\Theta\; | \; U(\theta')> U(\theta) + \tilde\epsilon\}) \leq \tilde\alpha\,\lambda(\Theta) \}\, .
\end{eqnarray*}
We have
$$
U(\theta') > U(\theta)  + \tilde\epsilon\quad\Leftrightarrow \quad
\rho\,U_\delta(\theta') > \rho\, [U_\delta(\theta) + \tilde\epsilon]
\quad\Rightarrow\quad \rho\, U_\delta(\theta')>U_\delta(\theta)
$$
which is proven by noticing that
$$
\rho\,[U_\delta(\theta) + \tilde\epsilon]\geq U_\delta(\theta)\,\,
\Leftrightarrow\,\,  (1-\rho) \geq U(\theta) (1-\rho)
$$
and $U(\theta)\in[0,\,1]$. Hence,
$$
 \{\theta'\in\Theta\; | \; \rho\, U_\delta(\theta')> U_\delta(\theta) \}
\quad \supseteq \quad
\{\theta'\in\Theta\; | \; U(\theta')> U(\theta) + \tilde\epsilon \} \, .
$$
Therefore,
$$
\pi_\delta (\{\theta'\in\Theta\; | \; \rho\, U_\delta(\theta')> U_\delta(\theta) \}) \leq {\bar\alpha}\quad\Rightarrow\quad
\pi_\delta (\{\theta'\in\Theta\; | \; U(\theta')> U(\theta) + \tilde\epsilon \}) \leq {\bar\alpha}\, .
$$
Let $Q_{\theta,\tilde\epsilon} :=  \{\theta'\in\Theta\; | \; U(\theta')> U(\theta) + \tilde\epsilon \}$
and notice that
$$
\pi_\delta(\{\theta'\in\Theta\; | \; U(\theta')> U(\theta) + \tilde\epsilon \})\,  =\,
\frac{
\mbox{$\displaystyle \int_{Q_{\theta,\tilde\epsilon}}
U(\theta')\lambda(d\theta') +\delta
\lambda(Q_{\theta,\tilde\epsilon})$}
}
{
\mbox{$\displaystyle \int_{\Theta} U(\theta')\lambda(d\theta') + \delta\lambda(\Theta)$}
}\, .
$$
We obtain
\begin{align*}
\pi_\delta(\{\theta'\in\Theta\; | \; U(\theta')> U(\theta) + \tilde\epsilon \}) \leq {\bar\alpha}
\, &\Rightarrow\,
\tilde\epsilon\,\lambda(Q_{\theta,\tilde\epsilon}) + \delta\lambda(Q_{\theta,\tilde\epsilon})
\leq {\bar\alpha} (1+\delta)\lambda(\Theta)  \\
&\Rightarrow\,\lambda
(\{\theta'\in\Theta\; | \; U(\theta')> U(\theta) +\tilde\epsilon \})  \leq \tilde\alpha\, \lambda(\Theta)\, .
\end{align*}
Hence we can conclude that
$$
\pi_\delta (\{\theta'\in\Theta\; | \; \rho\, U_\delta(\theta')> U_\delta(\theta) \})\leq {\bar\alpha}
\,\,\Rightarrow\,\,
\lambda (\{\theta'\in\Theta\; | \; U(\theta')> U(\theta) + \tilde\epsilon \}) \leq \tilde\alpha\, \lambda(\Theta)
$$
 and the second part of the proof is complete.

We have shown that given ${\bar\alpha}\in (0,\, 1]$, $\rho\in (0,\, 1]$,
$\tilde\epsilon := (\rho^{-1}-1)(1+\delta)$
and $\tilde\alpha := \frac{1+\delta}{\tilde\epsilon +\delta}\,{\bar\alpha}$, then
$$
\pi\left(\Theta(\tilde\epsilon,\tilde \alpha);J,\delta\right)\geq
\frac{1}
{\mbox{$\displaystyle 1 + \rho^{\,J}\frac{1-\bar\alpha}{\bar\alpha}\frac{1+\delta}{\delta}$}}
=\frac{1}
{\mbox{$\displaystyle
1+\left[\frac{1+\delta}{\tilde\epsilon+1+\delta}\right]^J
\left[\frac{1}{\tilde\alpha}\frac{1+\delta}{\tilde\epsilon+\delta}-1\right]
\frac{1+\delta}{\delta}$}}\, .
$$
Notice that $\tilde\epsilon\in [0,\, 1]$ and $\tilde\alpha\in(0,\, 1]$ are
linked through a bijective relation to
$\rho\in[\frac{1+\delta}{2+\delta},\, 1]$  and
$\bar\alpha\in (0,\,\frac{\tilde\epsilon+\delta}{1+\delta}]$.
Hence, the statement of the
theorem is eventually obtained by setting the desired
$\tilde\epsilon=\epsilon$ and $\tilde\alpha=\alpha$
in the above inequality.
\end{proof_ref}

To prove Corollary~\ref{cor:Jbound} we will need the following fact.

\begin{proposition}\label{prop:trivialfact}
For all $x>0$, $y>1$,
\[
\log\left(\frac{x+y}{y}\right)\geq
\frac{x}{x+y}\, .
\]
\end{proposition}
\begin{proof_text}
Fix an arbitrary $y>1$. If $x=0$ then
$\displaystyle
\log\frac{x+y}{y}=0=
\frac{x}{x+y}$.
Moreover,
\[
\frac{d}{dx}\left(\frac{x}{x+y}\right)
=\frac{y}{(x+y)^2}
\leq \frac{1}{x+y}
=\frac{d}{dx}
\left(\log\frac{x+y}{y}\right).\]
\end{proof_text}

\begin{proof_ref}{Corollary~\ref{cor:Jbound}}
To make sure that
$\pi(\Theta(\epsilon,\alpha);J,\delta) \geq \sigma$
we need to select $J$ such that
\[
\frac{1}{1+\left[\frac{1+\delta}{\epsilon+1+\delta}\right]^J
\left[\frac{1}{\alpha}\frac{1+\delta}{\epsilon+\delta}-1\right]
\frac{1+\delta}{\delta}}\geq \sigma,
\]
or, in other words,
\[
\left[\frac{\epsilon+1+\delta}{1+\delta}\right]^J \geq
\frac{\sigma}{1-\sigma}
\left[\frac{1}{\alpha}\frac{1+\delta}{\epsilon+\delta}-1\right]
\frac{1+\delta}{\delta}.
\]
It therefore suffices to choose $J$ such that
\[
\left[\frac{\epsilon+1+\delta}{1+\delta}\right]^J \geq
\frac{\sigma}{1-\sigma}
\frac{1}{\alpha}\left[\frac{1+\delta}{\delta}\right]^2.
\]
Taking logarithms
\[
J\log\frac{\epsilon+1+\delta}{1+\delta} \geq
\log\frac{\sigma}{1-\sigma}+
\log\frac{1}{\alpha}+
2\log\frac{1+\delta}{\delta}.
\]
Using the result of Proposition~\ref{prop:trivialfact}
with $x=\epsilon$ and $y=1+\delta$ one eventually obtains
that it suffices to select $J$ according to
inequality (\ref{eq:Jbound}).
\end{proof_ref}

\begin{proof_ref}{Proposition~\ref{prop:optidelta}}
Notice that
\begin{align*}
\frac{d}{d\delta}f(\delta)& =\frac{1}{\epsilon}
\left[\log\frac{\sigma}{1-\sigma}+\log\frac{1}{\alpha}+
  2\log\frac{1+\delta}{\delta}-
  2\frac{1+\epsilon+\delta}{\delta(1+\delta)}\right]\\
\frac{d^2}{d\delta^2}f(\delta)& =
2\frac{1+\epsilon+\delta+2\epsilon\delta}{\epsilon\delta^2(1+\delta)^2}>0
\end{align*}
and therefore the function $f(\delta)$ is convex in $\delta$.
The second equation ensures that if $f(\delta)$ attains a minimum then
it is unique. To complete the proof we need to show that the equation
\begin{equation}\label{eq:deriv}
\frac{d}{d\delta}f(\delta) =\frac{1}{\epsilon}
\left[\log\frac{\sigma}{1-\sigma}+\log\frac{1}{\alpha}+
  2\log\frac{1+\delta}{\delta}-
  2\frac{1+\epsilon+\delta}{\delta(1+\delta)}\right]=0
\end{equation}
always has a solution for $\delta>0$. To simplify the notation define
\begin{align*}
f_1(\delta) & = \log\frac{1+\delta}{\delta} +
\log\frac{\sqrt\sigma}{\sqrt{1-\sigma}}+\log\frac{1}{\sqrt{\alpha}}\\
f_2(\delta) & = \frac{1+\epsilon+\delta}{\delta(1+\delta)}.
\end{align*}
Then~\eqref{eq:deriv} simplifies to $f_1(\delta) = f_2(\delta)$.
It is easy to see that both $f_1$ and $f_2$ are monotone decreasing
functions of $\delta$ and
\begin{align*}
\lim_{\delta\rightarrow 0} f_1(\delta) & = \lim_{\delta\rightarrow 0}
f_2(\delta) = \infty,\\
\lim_{\delta\rightarrow \infty} f_1(\delta) & =
\log\frac{\sqrt\sigma}{\sqrt{1-\sigma}}+\log\frac{1}{\sqrt{\alpha}} >0
\mbox{ for $\sigma\in(0.5,\,1)$, }
\mbox{ and } \lim_{\delta\rightarrow \infty} f_2(\delta) = 0.
\end{align*}
Moreover, as $\delta$ tends to $0$, $f_1(\delta)$ tends to infinity
more slowly than $f_2(\delta)$ ($O(\log(1/\delta)$ instead of
$O(1/\delta)$). Therefore the two function have to cross for some
$\delta>0$.
\end{proof_ref}

\begin{proof_ref}{Proposition~\ref{prop:MTJb}}
Let $p(\cdot;J,\delta)$ denote the density of $\pi(\cdot;J,\delta)$.
Consider any $\theta^* \in \Theta^*$.
We have:
\begin{align*}
p(\theta;J,\delta) & =
\frac{\left[U(\theta)+\delta\right]^J}{\int_{\theta' \in \Theta}
  \left[U(\theta')+\delta\right]^J \lambda(d\theta')}\\
& \leq \frac{\left[U(\theta^*)+\delta\right]^J}{\int_{\theta' \in
    \Theta}\left[U(\theta')+\delta\right]^J \lambda(d\theta')}\\
& = \frac{\int_{\theta' \in \Theta^*}\left[U(\theta')+\delta\right]^J
  \lambda(d\theta')}{\lambda(\Theta^*)}
  \frac{1}{\int_{\theta' \in
    \Theta}\left[U(\theta')+\delta\right]^J \lambda(d\theta')}\\
& \leq \frac{\int_{\theta' \in \Theta}\left[U(\theta')+\delta\right]^J
  \lambda(d\theta')}{\lambda(\Theta^*)}
  \frac{1}{\int_{\theta' \in
    \Theta}\left[U(\theta')+\delta\right]^J \lambda(d\theta')}\\
& = \frac{\lambda(\Theta)}{\lambda(\Theta^*)}\frac{1}{\lambda(\Theta)}\\
& \leq \frac{1}{\beta}\frac{1}{\lambda(\Theta)}.
\end{align*}
Recall that the independent uniform proposal distribution over $\Theta$
has density $q(\theta)=\frac{1}{\lambda(\Theta)}$.
Hence, from the above inequality we obtain that
$M=\frac{1}{\beta}$ satisfies the inequality
in the statement of Theorem \ref{thm:MT}.
Therefore, we can write
$\|P_{\bm{\theta}_k} - \pi(\,\cdot\,\,;J,\delta)\|_{\mbox{\tiny TV}}\leq (1 -\beta)^k$.
Hence, $(1 -\beta)^k\leq\gamma\rho \Rightarrow
P_{\bm{\theta}_k}(\Theta(\epsilon,\alpha);J,\delta)\geq \rho$,
from which (\ref{eq:MTJb}) is eventually obtained.
 \end{proof_ref}

To prove Proposition \ref{prop:MTJ} we first establish a general fact.

\begin{proposition}\label{prop:JhatJ}

Let the notation and assumptions of Theorem \ref{thm:bound} hold.
Let $p(\cdot;J,\delta)$ denote the density of $\pi(\cdot;J,\delta)$.
For all $J \geq \hat{J} \geq 0$ and $\delta>0$
$$
p(\theta;J,\delta) \leq
\left(\frac{1+\delta}{\delta}\right)^{J-\hat{J}}
p(\theta;\hat{J},\delta),\quad\forall\theta \in \Theta.
$$
\end{proposition}
\begin{proof_text}
\begin{align*}
p(\theta;J,\delta) & =
\frac{\left[U(\theta)+\delta\right]^J}{\int_{\theta' \in \Theta}
  \left[U(\theta')+\delta\right]^J \lambda(d\theta')}\\
& =\frac{\left[U(\theta)+\delta\right]^J}
{\int_{\theta' \in \Theta}\left[U(\theta')+\delta\right]^J \lambda(d\theta')}
\frac{\int_{\theta' \in \Theta}\left[U(\theta')+\delta\right]^{\hat J} \lambda(d\theta')}
{\left[U(\theta)+\delta\right]^{\hat J}}
\frac{\left[U(\theta)+\delta\right]^{\hat J}}
{\int_{\theta' \in \Theta}\left[U(\theta')+\delta\right]^{\hat J} \lambda(d\theta')}\\
& = \left[U(\theta)+\delta\right]^{J-\hat{J}}
\frac{\int_{\theta' \in \Theta}
  \left[U(\theta')+\delta\right]^{\hat J}
  \lambda(d\theta')}{\int_{\theta' \in \Theta}
  \left[U(\theta')+\delta\right]^{J-\hat{J}}
  \left[U(\theta')+\delta\right]^{\hat J} \lambda(d\theta')}
\,p(\theta;\hat{J},\delta)\\
& \leq (1+\delta)^{J-\hat{J}}
\frac{\int_{\theta' \in \Theta} \left[U(\theta')+\delta\right]^{\hat J} \lambda(d\theta')}
{\delta^{J-\hat{J}}  \int_{\theta' \in \Theta}\left[U(\theta')+\delta\right]^{\hat J} \lambda(d\theta')}
\,p(\theta;\hat{J},\delta)\\
& = \left(\frac{1+\delta}{\delta}\right)^{J-\hat{J}}
p(\theta;\hat{J},\delta).
\end{align*}
\end{proof_text}

\begin{proof_ref}{Proposition~\ref{prop:MTJ}}
If we set $\hat J = 0$ in Proposition \ref{prop:JhatJ}
we obtain that
$$
M = \left(\frac{1+\delta}{\delta}\right)^J
$$
satisfies the inequality in the statement of Theorem~\ref{thm:MT}.
Hence, we obtain
$$
\|P_{\bm{\theta}_k} - \pi(\,\cdot\,\,;J,\delta)\|_{\mbox{\tiny TV}}\leq \left[1 - \left(\frac{\delta}{1+\delta}\right)^J \right]^k.
$$
Hence, it suffices to have
$$
\left(\frac{(1+\delta)^J-\delta^J}{(1+\delta)^J}\right)^k \leq \gamma\rho
$$
in order to guarantee  $\|P_{\bm{\theta}_k} - \pi(\,\cdot\,\,;J,\delta)\|_{\mbox{\tiny TV}}\leq \gamma\rho$.
Taking logarithms this becomes
$$
k \log\left(\frac{(1+\delta)^J}{(1+\delta)^J-\delta^J}\right) \geq \log\frac{1}{\gamma \rho}
$$
and, by applying Proposition~\ref{prop:trivialfact} with $x=\delta^J$ and $y=(1+\delta)^J-\delta^J$,
we eventually obtain
$$
k \geq \left(\frac{1+\delta}{\delta}\right)^J\log\left(\frac{1}{\gamma\rho} \right).
$$
Eventually, one obtains (\ref{eq:MTJ}) by
changing the base of the logarithms in the right-hand
side of (\ref{eq:Jbound}) from $e$ to $\frac{1+\delta}{\delta}$, and
by substituting $J$ with the so-obtained expression
in the right-hand side of the above inequality.
\end{proof_ref}

{
\bibliographystyle{unsrt}
\bibliography{reportbib}
}

\end{document}